\documentstyle[amssymb,12pt]{amsart}
\newtheorem{prop}{Proposition}[section]

\newtheorem{lemma}{Lemma}[section]
\newtheorem{thm}{Theorem}[section]
\newtheorem{corollary}{Corollary}[section]

\theoremstyle{remark}
\newtheorem{remark}{Remark}

\begin{document}
\newcommand{\nc}{\newcommand} \nc{\on}{\operatorname}
\nc{\pa}{\partial}
\nc{\cA}{{\cal A}}\nc{\cB}{{\cal B}}\nc{\cC}{{\cal C}}
\nc{\cE}{{\cal E}}\nc{\cG}{{\cal G}}\nc{\cH}{{\cal H}}
\nc{\cX}{{\cal X}}\nc{\cR}{{\cal R}}\nc{\cL}{{\cal L}}
\nc{\cK}{{\cal K}}\nc{\cM}{{\cal M}}\nc{\cT}{{\cal T}}
\nc{\sh}{\on{sh}}\nc{\Id}{\on{Id}}\nc{\Diff}{\on{Diff}}\nc{\Grass}{\on{Grass}}
\nc{\ad}{\on{ad}}\nc{\Der}{\on{Der}}\nc{\End}{\on{End}}\nc{\res}{\on{res}}
\nc{\Imm}{\on{Im}}\nc{\limm}{\on{lim}}\nc{\Ad}{\on{Ad}}\nc{\Den}{\on{Den}}
\nc{\Hol}{\on{Hol}}\nc{\Det}{\on{Det}}\nc{\Bun}{\on{Bun}}\nc{\tr}{\on{tr}}
\nc{\Cas}{\on{Cas}}\nc{\Hitch}{\on{Hitch}}\nc{\Ker}{\on{Ker}}
\nc{\id}{\on{id}}\nc{\Aut}{\on{Aut}}\nc{\Hom}{\on{Hom}}
\nc{\diag}{\on{diag}}\nc{\corr}{\on{corr}}\nc{\Fil}{\on{Fil}}
\nc{\de}{\delt}\nc{\si}{\sigma}\nc{\ve}{\varepsilon}
\nc{\al}{\alpha}
\nc{\CC}{{\mathbb C}}\nc{\ZZ}{{\mathbb Z}}\nc{\NN}{{\mathbb N}}
\nc{\PP}{{\mathbb P}}\nc{\VV}{{\mathbb V}}\nc{\WW}{{\mathbb W}}
\nc{\FF}{{\mathbb F}}
\nc{\AAA}{{\mathbb A}}\nc{\cO}{{\cal O}} \nc{\cF}{{\cal F}}\nc{\cD}{{\cal D}}
\nc{\la}{{\lambda}}\nc{\G}{{\mathfrak g}}\nc{\A}{{\mathfrak a}}
\nc{\HH}{{\mathfrak h}}\nc{\dd}{{\mathfrak d}}\nc{\Heis}{{\mathfrak H}}
\nc{\N}{{\mathfrak n}}\nc{\B}{{\mathfrak b}}
\nc{\La}{\Lambda}\nc{\ul}{{\underline{\ell}}}
\nc{\uQ}{{\underline{Q}}}\nc{\uchi}{{\underline{\chi}}}
\nc{\ula}{{\underline{\la}}}\nc{\uP}{{\underline{P}}}\nc{\uV}{{\underline{V}}}
\nc{\g}{\gamma}\nc{\eps}{\epsilon}\nc{\wt}{\widetilde}
\nc{\wh}{\widehat}
\nc{\bn}{\begin{equation}}\nc{\en}{\end{equation}}
\nc{\SL}{{\mathfrak{sl}}}\nc{\GL}{{\mathfrak{gl}}}
\nc{\SG}{{\mathfrak{S}}}

%
%
%

\newcommand{\ldar}[1]{\begin{picture}(10,50)(-5,-25)
\put(0,25){\vector(0,-1){50}}
\put(5,0){\mbox{$#1$}} 
\end{picture}}

\newcommand{\lrar}[1]{\begin{picture}(50,10)(-25,-5)
\put(-25,0){\vector(1,0){50}}
\put(0,5){\makebox(0,0)[b]{\mbox{$#1$}}}
\end{picture}}

\newcommand{\luar}[1]{\begin{picture}(10,50)(-5,-25)
\put(0,-25){\vector(0,1){50}}
\put(5,0){\mbox{$#1$}}
\end{picture}}

\title[Hecke-Tyurin parametrization of the Hitchin and KZB systems]
{Hecke-Tyurin parametrization of the Hitchin and KZB systems}

\author{B. Enriquez}

\address{B.E.: DMA (UMR 8553 du CNRS), ENS, 45 rue d'Ulm, 75005 Paris, France}

\author{V. Rubtsov}

\address{V.R.: D\'ept.\  de Math\'ematiques, Universit\'e d'Angers, 
Angers, France}

\address{ITEP, 25, Bol. Cheremushkinskaya, Moscow 117259, Russia}

\date{November 1999}

\begin{abstract}

We study the parametrization of the moduli space $\Bun_2(C)_\cL$  of
rank $2$ bundles over a curve $C$ with fixed determinant, provided  by
Hecke modifications at fixed points of the trivial bundle. This
parametrization is closely related to the  Tyurin parametrization of
vector bundles over curves.  We use it to parametrize the Hitchin and
KZB systems, as well as lifts of the Beilinson-Drinfeld $\cD$-modules. 
We express a generating series  for the lifts of the Beilinson-Drinfeld 
operators in terms of a ``quantum $L$-operator'' $\ell(z)$.    We
explain the relation to earlier joint work with G. Felder, based on
parametrization by flags of  bundles, and introduce filtrations on
conformal blocks, related with Hecke  modifications.      

\end{abstract}

\maketitle

\section{Introduction}

\subsection{Tyurin maps and Hecke parametrizations}

In \cite{Tyurin}, A.\ Tyurin introduced the following parametrization of
the moduli space $\Bun_n(C)_{ng}$ of vector bundles of rank $n$ and
of degree $ng$ over a curve $C$ of genus  $g = g(C) \geq 1$. The space
$H^0(C,\cE)$ of sections of a generic such bundle $\cE$ is of dimension
$n$. For $P$ a point of $C$, denote by $\cE_P$ the fiber of $\cE$ at $P$ and  define 
$\Sigma_P$ as the subspace of $\cE_P$ generated by the sections of $\cE$. 
$\Sigma_P$ is equal to the fiber of  $\cE_P$ at $P$,
except when $P$ is one of $ng$ points of $C$, $P_1,\ldots,P_{ng}$.
Generically,  the $P_i$ are distinct and $\Sigma_P$ is of codimension
$1$ in  $\cE_P$. Let us denote by $\ell_i$ the subspace of 
$H^0(C,\cE)^*$ generated by evaluation at $P_i$. Define $\cT$ as the map
 assigning to $\cE$, the collection $(P_i)_{i = 1,\ldots,ng}$ together with the collection
of  lines $(\ell_i)_{i = 1,\dots,ng}$ in $H^0(C,\cE)^*$, modulo diagonal action of
$SL_n(\CC)$. The resulting map  
$$
\cT : \Bun_n(C)_{ng} \to [C \times \PP
(\CC^n)]^{(ng)} / SL_n(\CC)
$$
is called the Tyurin map (we denote by $V^{(N)}$ the $N$th symmetric power of a 
variety $V$). After its introduction in \cite{Tyurin}, the Tyurin map was used 
in works on integrable differential systems (see \cite{Grin,Krichever,PW}). 

The inverse of the map $\cT$ can be naturally described in terms of 
Hecke modifications. For $\cE$ a rank $n$ vector bundle over $C$, $\uP =
(P_i)_{i = 1,\ldots,N}$ a family of points of $C$ and $\uV = (V_i)_{i =
1,\ldots,N}$ a family of vector subspaces of
$\cE_{P_1},\ldots,\cE_{P_N}$,   define the Hecke modification
$\cH(\cE,\uP,\uV)$  of $\cE$ as the bundle whose  space of sections over
an open subset $U\subset C$ is the space  of rational sections of 
$\cE$ over $U$, regular except at the points $P_i$, where they may have simple
poles,  with residue in $V_i$. In the case where $\cE$ is the trivial
bundle $\cO_C^n$,  and $V_i$ are lines $\ell_i$ in $\CC^n$, let us set
$\ul = (\ell_i)_{i = 1,\ldots, N}$ and  $\cH(\uP,\ul)  =
\cH(\cO_C^n,\uP,\ul)$.  Since the automorphism group of $\cO_C^n$ is
equal to $GL_n(\CC)$, the map 
$$
\cH : (C \times \PP(\CC^n))^{(N)} \to \Bun_n(C)
$$ 
factors through the diagonal action of  $GL_n(\CC)$ (or
$SL_n(\CC)$). In the case where $N = ng$, the resulting map  (also
denoted $\cH$)  is inverse to the Tyurin map. 

In the case $n = 1$,  the map $\cH$ coincides with the Abel-Jacobi
map  from the symmetric powers of $C$ to its Picard variety.  It is
therefore natural to view $\cH$ as a nonabelian analogue of the
Abel-Jacobi map. 

\subsection{Main results}

In this paper, we will keep the points $\uP = (P_i)_{i = 1,\ldots,ng}$
fixed, and associate to them the map $\cH_\uP$ such that $\cH_\uP(\ul) =
\cH(\uP,\ul)$.  We will study $\cH_\uP$ (sect.\ \ref{sect:1}) and 
parametrize in terms of it natural objects attached to  the moduli
spaces of vector bundles, in the case $n=2$. These objects are the 
Hitchin system (sect.\ \ref{sect:hitchin}), the bundle of conformal
blocks (sect.\ \ref{sect:cb}),  the Beilinson-Drinfeld $\cD$-modules
(sect.\ \ref{sect:bd}) and the Knizhnik-Zamolodchikov-Bernard (KZB) connection
(sect.\ \ref{sect:kzb}); they were respectively introduced in 
\cite{Hitchin}, \cite{KZ,Bernard,TUY} and \cite{BD}. The main results of
this paper are the following. 

\subsubsection{Properties of $\cH_\uP$}

Let us set $N = 3g, G = SL_2$. $\cH_\uP$ is then a map from $(\CC P^1)^{3g} / 
G(\CC)$ from the moduli space $\Bun_2(C)_{\cO(\sum_i P_i)}$
of rank $2$ bundles with determinant isomorphic to $\cO(\sum_i P_i)$. 

Let us fix a choice of $a$- and $b$-cycles $(A_a,B_a)_{a = 1,\ldots,g}$
on $C$, and let $(\omega_a)_{a = 1,\ldots, g}$ be the associated 
basis of regular differentials on $C$; we have $\int_{A_a} \omega_b
= 2i\pi \delta_{ab}$.

Let us set 
\begin{align} \label{expr:Den}
 \Den(\uP,\ul) = & \sum_{\sigma\in \SG_{3g}: 
 \sigma(3a -2) < \sigma(3a -1) < \sigma(3a) \on{\ for\ }
 a = 1,\ldots,g}
 \\ & \nonumber 
 \varepsilon(\sigma) \sigma \cdot \{ \prod_{a = 1}^g
 \omega_a(P_{3a -2})\omega_a(P_{3a -1})\omega_a(P_{3a})
 \ell_{3a -2;3a -1}\ell_{3a -2;3a}\ell_{3a -1;3a}\} , 
\end{align}
where $\ell_{ij} = \ell_i - \ell_j$, and the elements $\sigma$ of 
$\SG_{3a}$ act by permuting both the $P_i$ and the $\ell_i$.

\begin{thm} \label{thm:local:isom} The preimage of 
$\Bun_2(C)_{\cO(\sum_i P_i)}^{stable}$ by $\cH_\uP$ contains the open  
$\{\ul\in (\CC P^1)^{3g}|\Den(\uP,\ul) \neq 0\} /
G(\CC)$. Moreover, the  restriction of  $\cH_{\uP}$ to
$$
\{\ell\in (\CC P^1)^{3g}|\Den(\uP,\ul) \neq 0\} / G(\CC)
$$  
is an etale morphism from this variety  to its image, which is a
dense subset of  $\Bun_2(C)_{\cO(\sum_i P_i)}^{stable}$.  
\end{thm}

\subsubsection{Parametrization of the Hitchin system}

By symplectic reduction,  we may identify $T^*[(\CC P^1)^{3g}]/G$ with
the quotient of  $\{(\ul,\ula)\in \CC^{3g}\times\CC^{3g} | \sum_i \la_i =  
\sum_i \ell_i \la_i = \sum_i \ell_i^2\la_i = 0 \} / G$ (see
sect.\ \ref{sect:hitchin}). Let $(e,f,h)$ be the Chevalley basis of 
$\bar\G = \SL_2$ and for $(\ul,\ula)$ in $\CC^{3g} \times \CC^{3g}$, 
let us set 
\begin{align} \label{def:A}
& \nonumber A(\ul,\ula|z) =  \sum_{i = 1}^{3g}
\la_i(-e + \ell_i h + \ell_i^2 f) \omega^{(P_i)}(z) 
+ {1\over{\Den(\uP,\ul)}} 
\sum_{i,j = 1,\ldots,3g;j\neq i} \la_i \ell_{ij}^2 \omega^{(P_i)}(P_j)
\\ & 
\{  {1\over 2} \pa^2_{\ell_j}\Den(\uP,\ul)_{| \ell_j = 0} e
+ {1\over 2} \pa_{\ell_j}\Den(\uP,\ul)_{| \ell_j = 0} h  
- \Den(\uP,\ul)_{| \ell_j = 0} f \}_{|P_j = z} , 
\end{align}
where $\ell_{ij} = \ell_i - \ell_j$ and 
$$
\omega^{(z)}(P) = d_P \ln\Theta( A(P) - A(z) + \delta - \Delta) 
$$
is the Green function on $C$ (see the conventions for theta-functions in 
sect.\ \ref{sect:theta}). 
$A(\ul,\ula|z)$ is an element of $\bar\G\otimes\Omega(C)$, where
$\Omega(C)$ is the space of rational differentials on $C$.  
Let us set $H(\ul,\ula|z) = \tr\rho[A(\ul,\ula|z)]^2$. 

The Hitchin Hamiltonians $\Hitch_\al$ are regular maps defined on  $T^*
\Bun_2(C)_{\cO(\sum_i P_i)}$ (see sect.\ \ref{sect:hitchin} and
\cite{Hitchin}).

\begin{prop} \label{prop:hitchin}
$z\mapsto H(\ul,\ula|z)$ is a regular quadratic differential
on $C$. Let us set 
$$
H(\ul,\ula|z) = \sum_{\al = 1}^{3g - 3} H_\al(\ul,\ula)
\omega^{(2)}_\al(z). 
$$
The $H_\al(\ul,\ula)$ are rational functions in $\ell_i$ and
quadratic polynomials in $\la_i$. They are $G$-invariant functions on 
$\{(\ul,\ula) =  (\ell_i,\la_i)_{i = 1,\ldots, 3g} | 
\sum_i \la_i = \sum_i \la_i \ell_i = \sum_i \la_i
\ell_i^2 = 0\}$, and define therefore a family of functions on 
$T^*[(\CC P^1)^{3g} / G(\CC)]$. These functions form a Poisson-commutative
family  coinciding with the $\Hitch_\al\circ T^*\cH_{\uP}$. 
\end{prop}

\subsubsection{Explicit form of the Beilinson-Drinfeld operators}

The Beilinson-Drinfeld operators are a commuting family of globally 
defined differential operators 
$(T^{BD}_{\al})_{\al = 1,\ldots,3g-3}$ 
on $\Bun_2(C)_{\cO(\sum_i P_i)}$, 
twisted by $\det^{-2}$, where $\det$ is the determinant line bundle 
over $\Bun_2(C)_{\cO(\sum_i P_i)}$ (see \cite{BD} and sect.\ \ref{sect:bd}). 

\begin{prop} \label{prop:image:det}
 $\cH_\uP^*(\det)$ is isomorphic with $\cO(1)^{\boxtimes 3g}$. 
\end{prop}

Let $a_\uP(\ul,z)$ be the element of $\bar\G\otimes 
\CC[\ell_i,\pa_{\ell_i}, i = 1,\ldots,3g]
[\Den(\uP,\ul)^{-1}]\otimes \Omega(C) $ equal to 
\begin{align*}
& a_\uP(\ul,z) = - {1\over{\Den(\uP,\ul)}}\sum_i G(P_i,z) dP_i 
\{ ({1\over 2} h + \ell_i f ) \Den(\uP,\ul)_{|\ell_i = 0} 
\\ &  
- (e + \ell_i^2 f) \pa_{\ell_i}\Den(\uP,\ul)_{|\ell_i = 0} 
+ ( 2 \ell_i e - \ell_i^2 h ) \pa^2_{\ell_i}\Den(\uP,\ul)_{|\ell_i = 0} 
\}_{P_i = z} , 
\end{align*}
and let $s_\uP(\ul,z)$ be the element of 
$\CC[\ell_i,\pa_{\ell_i},i = 1,\ldots,3g]
[\Den(\uP,\ul)^{-1}] \otimes
\Omega^2(C)$ ($\Omega^2(C)$ is the space of rational sections of $\Omega_C^{\otimes 2}$) 
given by 
\begin{align*}
s_\uP(\ul,z) = & {k\over{\Den(\uP,\ul)}}
\sum_{i = 1}^{3g} d_z[G(P_i,z)dP_i]
\\ & \{ - {1\over 2} \Den(\uP,\ul)_{|\ell_i = 0}
+ \ell_i \pa_{\ell_i} \Den(\uP,\ul)_{|\ell_i = 0}
- \ell_i^2 \pa^2_{\ell_i} \Den(\uP,\ul)_{|\ell_i = 0} \}_{P_i = z}.  
\end{align*}

Let us also set $G(z,w)dz = \omega^{(z)}(w)$,  
\begin{align*}
& \mu_{i,\uP}(\ul,z) =  (e - \ell_i h - \ell_i^2 f) G(z,P_i) dz  
 + {1\over{\Den(\uP,\ul)}} \cdot \\ & 
 \cdot \sum_{j\neq i} \ell_{ij}^2
[ \pa_{\ell_j}^2\Den(\uP,\ul)_{|\ell_j = 0} e
- {1\over 2} \pa_{\ell_j}\Den(\uP,\ul)_{|\ell_j = 0} h
- {1\over 2} \Den(\uP,\ul)_{|\ell_j = 0} f ]_{P_j = z}
G(P_j,P_i) dP_j , 
\end{align*} 
\begin{align*}
&  \nu_\uP(\ul,z) = - \sum_i ({1\over 2} h + \ell_i f) G(z,P_i)dz 
+ {1\over{\Den(\uP,\ul)}}  
\sum_{i,j; i\neq j} G(P_j,P_i) dP_j   \ell_{ji} \cdot 
\\ & \cdot 
[ - \pa^2_{\ell_j}\Den(\uP,\ul)_{|\ell_j = 0} e  
+ {1\over 2} \pa_{\ell_j}\Den(\uP,\ul)_{|\ell_j = 0} h  
+ {1\over 2} \Den(\uP,\ul)_{|\ell_j = 0} f 
]_{P_j = z} 
\end{align*}
and 
$$
\ell^{diff}_\uP(z) := \sum_i \mu_{i,\uP}(\ul,z) \pa_{\ell_i} - k \nu_\uP(\ul,z) .   
$$
$\ell^{diff}_\uP(z)$ is an element of  $\bar\G\otimes
\CC[\ell_i,\pa_{\ell_i}, i = 1,\ldots,3g] [\Den(\uP,\ul)^{-1}]\otimes\Omega(C)$.

Let us set $\sum_\al x_\al \otimes y_\al = e\otimes f + f\otimes e 
+ {1\over 2} h\otimes h$ and  
$$ 
T^{diff}_\uP(z) = 
\sum_\al \ell^{diff}_{\uP,x_\al}(z)\ell^{diff}_{\uP,y_\al}(z) 
+ \sum_\al a_{\uP,x_\al}(\ul,z)\ell^{diff}_{\uP,y_\al}(z) + s_\uP(\ul,z) . 
$$

\begin{thm} \label{thm:bd}
Assume that $k+2 = 0$. Then $T^{diff}(z)$ is a regular quadratic differential 
on $C$. Let us $(\omega^{(2)}_\al(z))_{\al = 1,\ldots, 3g-3}$ be a basis of 
$H^0(C,\Omega_C^2)$ and set 
$$
T^{diff}_\uP(z) = \sum_{\al = 1}^{3g - 3} T^{diff}_{\uP,\al} \omega^{(2)}_\al(z) 
$$
The $T^{diff}_{\uP,\al}$ form a commuting family of 
differential  operators. 

For $\sigma$ a section of $\det^{-2}$ on an open  subset of the stable
locus, we have 
$$
T^{diff}_{\uP,\al} (\cH_\uP^* \sigma)  =  \cH_\uP^* [T^{BD}_{\al}
(\sigma)] .  
$$
\end{thm}

\subsubsection{Functional parametrization of conformal blocks}

Let us denote by $A : C \to J^1(C)$ the Abel-Jacobi map of $C$  (see
sect.\ \ref{AJ}).  Assume that we have $\sum_i A(P_i) = 3g A(P_0)$.
$\Bun_2(C)_{\cO(\sum_i P_i)}$ is then canonically isomorphic to 
$\Bun_2(C)_{\cO(3g P_0)}$. Let $z_{P_0}$ be a local coordinate at 
$P_0$, let $\cK = \CC((z_{P_0}))$ and $\cO = \CC[[z_{P_0}]]$. 
Then we construct an element 
$$
g_{\uP,\ul} \in GL_2(\cK),
$$ 
such that its class in the double  quotient 
$GL_2(H^0(C \setminus \{P_0\}, \cO_C)) 
\setminus GL_2(\cK) / GL_2(\cO)
 = \coprod_{l\in\ZZ} \Bun_2(C)_{\cO(lP_0)}$ corresponds to that of 
$\cH_\uP(\ul)$ (Lemma \ref{lemma:pties:gPl}). 

Let us set $G = SL_2$. Let $\wt T_{G(\cK)}$ be the automorphism of 
$G(\cK)$ defined as conjugation by $\pmatrix z_{P_0} & 0 \\ 0 & 1
\endpmatrix$. Let $G(\cK) \rtimes \ZZ$ be the semidirect 
product of $G(\cK)$ by this automorphism. Then we have a sequence
of group inclusions
$$
G(\cK) \subset G(\cK)\rtimes \ZZ \subset GL_2(\cK) , 
$$
and the analogous sequence for centrally extended groups $\wh{G(\cK)}
\subset \wh{G(\cK)}\rtimes \ZZ \subset \wh{GL_2(\cK)}$,  (sect.\
\ref{semidirect}). (The advantage of working with $G(\cK)\rtimes \ZZ$
rather than $GL_2(\cK)$ is that is does not  require the introduction of
an additional Heisenberg algebra.) Moreover, $GL_2(\cK)$ is the product
of  $G(\cK)\rtimes \ZZ$ by its diagonal subgroup. We construct then  an
element $\wt g_{\ul,\uP}$ of  $G(\cK)\rtimes \ZZ$, which is the product
of $g_{\ul,\uP}$ by a  diagonal element of $GL_2(\cK)$.

In sect.\ \ref{rep:semidirect}, we also classify the integrable 
irreducible representations of $G(\cK)\rtimes \ZZ$. Let $\WW$ 
be such a representation. $\WW$ can be constructed as the quotient of 
a module induced by a $\bar\G$-module $W$; we have then an inclusion 
$W\subset \WW$. The space of conformal blocks of $\WW$ is then the 
space $CB(\WW)$ of linear forms on $\WW$, invariant by 
$\bar\G\otimes H^0(C\setminus\{P_0\},\cO_C)$.

Let $\psi$ be an element of $CB(\WW)$. We associate to it the 
quantity 
$$
f_\psi(\ul,\uP|v) = \langle \psi, \wt g_{\ul,\uP} v\rangle,   
$$
which we view as a function from $(\CC P^1)^{3g}$ to $W^*$.

\begin{prop} \label{prop:fun}
(functional properties of $f_\psi(\ul,\uP|v)$) For any $\uP$, the map
$v\mapsto f_\psi(\ul,\uP|v)$ belongs to $\Hom_{\bar\G} ( W, \CC[\ell_1]_{\leq k} 
\otimes \cdots \otimes \CC[\ell_{3g}]_{\leq k})$. 
\end{prop}

Here $\CC[\ell]_{\leq k}$ is the space of polynomials in $\ell$ of
degree $\leq k$; it is endowed with its standard $\bar\G$-module 
structure.  

As we will see, the map $\psi\mapsto f_\psi$ should be considered as a
``quantization'' of the cotangent to the Hecke map $T^*\cH_\uP$. 

The variation of $f_\psi(\ul,\uP|v)$ with $\uP$ is described as follows: 

\begin{prop} \label{variation:f:psi} 
Let $\wt f_\psi(\uP,\ul|v)$ be
a function on $C^{3g}\times (\CC P^1)^{3g}$ with values in $W^*$, whose restriction to 
$(C^{3g})_{P_0}\times (\CC P^1)^{3g}$ coincides with $f_\psi(\uP,\ul|v)$. We have 
$$
\pa_{P_i} \wt f_\psi(\ul,\uP|v) = (\La_i  \wt f_\psi)(\ul,\uP|v) + \sum_a \omega_a(P_i) 
(\cA_a \wt f_\psi)(\ul,\uP|v), 
$$
where 
\begin{align*}
& \La_i =  \sum_{j : j\neq i} r^{(P_i)}(P_j) \ell_{ji} 
{{\ell_j}\over{\ell_i}} ({{\pa}\over{\pa\ell_j}} - {k\over{\ell_j}}) 
\\ & - {1\over{\Den(\uP,\ul)}} 
\sum_{j,l : j\neq l} r^{(P_l)}(P_j) \ell_{jl}^2
{1\over{\ell_i}}
\Den(\uP,\ul)_{|\ell_l = 0, P_l \to P_i}
({{\pa}\over{\pa\ell_j}} - {k\over{\ell_j}}) 
\\ & 
- k {1\over{\Den(\uP,\ul)}}\sum_{j,l : l\neq j}
{{\ell_j^2}\over{\ell_i \ell_l}} r^{(P_j)}(P_l)
\Den(\uP,\ul)_{| \ell_j = 0, P_j \to P_i}
\\ & + k  ({1\over {\ell_i}} e - {1\over 2}h)
+ {k\over{\Den(\uP,\ul)}} \sum_j - {{\ell_j}\over{\ell_i}}
\Den(\uP,\ul)_{|\ell_j = 0, P_i \to P_j} [{1\over{\ell_j}} e - h - \ell_j f] , 
\end{align*}
where we denote by $x\in \bar\G$ the operator $f(\uP,\ul |v)\mapsto f(\uP,\ul|
\rho_W(x)v)$, 
and $\cA_a$ are some linear operators.  
\end{prop}

Here the term $\sum_a \omega_a(P_i) (\cA_a \wt f_\psi)(\ul,\uP|v)$ is the
Lagrange multiplier corresponding to the relation $\delta(\sum A(P_i)) =
0$. 

\subsubsection{Explicit form of the KZB connection}

Let us denote by $\cM_{g,1^\infty}$ the moduli space of triples
$m_\infty = (C,P_0,t)$ of a  smooth curve $C$ of genus $g$,  a marked
point $P_0$ of $C$ and a formal coordinate $t$ at $P_0$.  

$\cM_{g,1^\infty}$ may be identified with the set of all inclusions
$R\subset \CC((t))$, where $R$ is the coordinate ring of a smooth curve of 
genus $g$ minus a point. 

Let us set $\G = \bar\G\otimes\CC((t)) \oplus\CC K$. Consider $\WW$ as a 
$\G$-module and $\G^{out,ext}$ as a Lie subalgebra of $\G$. 

The bundle of conformal blocks on $\cM_{g,1^\infty}$ is defined as the subbundle of 
the constant bundle over $\cM_{g,1^\infty}$ with fiber $\WW^*$, consisting of the
$\G^{out,ext}$-invariant forms. We denote this bundle by $CB(\WW)$. 

A flat connection is defined on $CB(\WW)$ as follows (see \cite{TUY,Bernard,EF:sols}). 
Recall first that the map 
$$
\CC((t)) {{\pa}\over{\pa t}} \to \on{Vect}(\cM_{g,1^\infty}), \quad \xi\mapsto [\xi], 
$$
defined by $\exp(\eps[\xi])([R\subset \CC((t))]) = [(1 +
\eps\xi)(R)\subset \CC((t))]$, when $\eps^2 = 0$, induces a Lie algebra
morphism from  vector fields on the formal disc to vector fields on
$\cM_{g,1^\infty}$.  For $m_\infty\mapsto \psi(m_\infty)$ a local
section of $CB(\WW)$, let us set 
$$
\nabla_{[\xi]}^{CB}\psi(m_\infty) = \pa_{[\xi]} \psi(m_\infty) 
- \psi(m_\infty) \circ T[\xi] , 
$$
where $T[\xi] = \res_{z = P_0}(\xi(z)T(z))$.

Let us denote by $\cF$ the bundle over $\cM_{g,1^\infty}$, whose fiber $\cF(m_\infty)$ 
at $m_\infty$ is the space of functions $f(\uP,\ul|v)$ from $(C^{3g})_{P_0} \times 
(\CC P^1)^{3g}$ to $W^*$, such that 
$$
\pa_{P_i}f(\uP,\ul|v) = \La_i f(\uP,\ul|v) + \sum_a \omega_a(P_i) f_a(\uP,\ul|v), 
$$
where $f_a(\uP,\ul|v)$ are  functions  from $(C^{3g})_{P_0} \times (\CC
P^1)^{3g}$ to $W^*$. 

Define $\cM_{g,1^\infty,(3g)}$ as the moduli space of all pairs $(m_\infty,\uP)$,
such that $\sum_i A(P_i) = 3g A(P_0)$. Let $\xi$ be a formal vector field $\xi\in \CC((t)) 
{\pa\over{\pa t}}$, such that 
$$
\res_{P_0} (\xi\omega_a {{db_\uP}\over{b_\uP}}) = 0,  
$$ 
for each $a = 1,\ldots, g$. Then $\xi$ induces a vector field $[\xi]$ at $(m_\infty,\uP)$, 
defined by $(1 + \epsilon [\xi] )(m_\infty,\uP) = (m'_\infty,\uP')$, where 
$m'_\infty$ in the inclusion $(1 + \eps\xi)(R)\subset \CC((t))$, and $P'_i$ correspond the 
characters $\chi_{P_i}\circ (1 - \eps\xi)$ of $R$, where $\chi_{P_i}$ is the character of $R$
corresponding to $P_i$  
(see \cite{EF:sols}, Lemma 5.1). 

\begin{prop} \label{prop:kzb}
The rule 
$$
(\nabla^{\cF}_{[\xi]} f)(m_\infty,\uP,\ul|v) = 
\pa_{[\xi]} f(m_\infty,\uP,\ul|v) - \res_{z = P_0} (\xi(z) 
T^{diff}(z) f(m_\infty,\uP,\ul|v))  
$$
where $f(m_\infty,\uP,\ul|v))$ is a local section of $\cF$, defines a 
flat connection  $\nabla^{\cF}$ on $\cF$. 

Define $\corr$ as the bundle morphism from $CB(\WW)$  to $\cF$, 
such that $\corr(\psi)(m_\infty,\uP,
\ul|v) = \langle \psi , \rho_\WW^{group}(\wt g_{\uP,\ul})(v)\rangle$. 
Then $\nabla^{\cF} \circ \corr = \corr \circ \nabla^{CB}$. 
\end{prop}

We end the paper with the construction of a filtration naturally 
attached to $\cH_\uP$ (sect.\ \ref{sect:filt}), and questions about its
possible  relation with formulas of  \cite{Feigin:Loktev,Leclerc}. 

\subsection{Generalization to arbitrary semisimple groups} 

The above parametrization maps are generalized as
follows.  Let $G$ be a reductive complex group. Let $\cO$ be a formal
series ring $\CC[[t]]$,  and let $\cK$ be its fraction field $\CC((t))$.
Let $\Grass$ be the affine  Grassmann variety  $\Grass =  G(\cK) /
G(\cO)$. $\Grass$  is the union of all $\Grass_\chi = G(\cO) w_\chi G(\cO)
/ G(\cO)$,  where $\chi$ belongs to the semigroup $P_+$ of dominant coweights
of the Langlands dual group $\null^L G$ of $G$, and $w_\chi$ is
associated translation of the affine Weyl group of $G$. Let $\uchi =
(\chi_i)_{i = 1,\ldots,N}$ be a family of  elements of $P_+$. We have then 
Hecke maps  
$$
\cH_{\uP,\uchi} : \prod_{i = 1}^N \Grass_{\chi_i} \to \Bun_G(C), 
$$
obtained from the natural map $G^{out,\uP} \setminus
\prod_{i = 1}^N \Grass_{\chi_i} \to \Bun_G(C)$, where 
$G^{out,\uP}$ is the group of regular maps from $C \setminus \{P_i, i = 1,\dots, N\}$
to $G$. 

Moreover, the map $[\on{class\  of\ }g w_{\chi} g'] \mapsto [\on{class\
of\ }g]$ defines  isomorphisms  
\begin{equation} \label{isoms}
\Grass_\chi \to G(\cO) / [\null^{w_\chi} G(\cO)\cap
G(\cO)]  = N(\cO) / [\null^{w_\chi} N(\cO)\cap N(\cO)],
\end{equation}
where $N$ is the positive unipotent subgroup of $G$, and we set 
$\null^g h = ghg^{-1}$.

\subsection{Comparison with earlier work}

In \cite{EF}, G.\ Felder and one of us introduced another 
parametrization of the KZB connection. Let us explain its relation to  
the approch of the present work. The work \cite{EF} was based on a
variant of the Feigin-Stoyanovski description of the conformal blocks, which
relies on the map 
$$
FS : \Bun_B(C)_{P_0} \to \Bun_G(C).
$$
Here $P_0$ is a fixed point of $C$ and $\Bun_B(C)_{P_0}$ is the  union
of all $\Bun_B(C)_{P_0,\chi}$, which is the moduli space of all 
$B$-bundles over $C$, with associated $T$-bundle isomorphic to
$\oplus_{i\ \on{simple}} \cO(\langle \chi, \al_i\rangle P_0)$.  Here $B$
the Borel subgroup  of $G$ containing $N$, and $T$ is the associated
Cartan subgroup. Let us denote by  $FS_\chi$ the restriction of $FS$ to 
$\Bun_B(C)_{P_0,\chi}$. 

The Hecke map $\cH_{P_0,\chi}$ factors through the above
map $FS_\chi$. 
Indeed, the composition of $\cH_{P_0,\chi}$ with the isomorphism 
(\ref{isoms}) is the map 
\begin{equation} \label{klez}
N(\cO) / [\null^{w_\chi} N(\cO)\cap N(\cO)] \to \Bun_G(C) , 
\end{equation}
class of $n\mapsto$ class of $n w_\chi$. On the other hand,   
$\Bun_B(C)_{P_0,\chi}$ is equal to $ N(R) \setminus \{N(\cK) w_\chi\} /
N(\cO)$, where $R = H^0(C\setminus \{P_0\},\cO_C)$. 
The map (\ref{klez}) therefore factors through  $\Bun_B(C)_{P_0,\chi}$. 

It is probable that the maps $\cH_{\uP,\uchi}$ have the following 
behavior. When $P_i$ and $P_j$ tend to some point $P$, the limit of  
$\cH_{\uP,\uchi}$ should be related to $\cH_{\uP',\uchi'}$, 
where $\uP' = ((P_\al)_{\al\neq i,j}, P)$ and $\uchi' = 
((\chi_\al)_{\al\neq i,j}, \chi_i + \chi_j)$. One could then think 
that the maps $FS$, which involve all coweights of $\null^L G$, 
can be obtained as limits of the maps $\cH_{\uP,\uchi}$, where 
$\chi_i$ are simple coweights of $\null^L G$, and that one can obtain
the differential systems expressed in \cite{EF} from those of 
the present work by a similar limiting procedure. 

\subsection{Open problems} 

1) It is easy to write variants of our parametrization of the
Beilinson-Drinfeld  $\cD$-modules, corresponding to the parametrization
$(\CC P^1)^{3g - 3} \to \Bun_2(C)_{\cO(\sum_i P_i + \sum_j Q_j)}$, 
$\ul\mapsto \cH_{\uP,\uQ}(\ul,\ul^Q)$, where $(Q_j)_{j = 1,\ldots, p}$
are fixed points of $C$ and $\ul^Q$ a fixed element of $(\CC P^1)^p$. 
Using these variants, one could hope to check explicitly the Hecke
eigenvalue  property of the Beilinson-Drinfeld $\cD$-modules.  

2) It should also be possible to write the analogues of our systems in the 
$\SL_n$ case; it that case the parametrization maps could be $(\CC P^n)^{n+1}
\to \Bun_n(C)_\cL$. 

3) In \cite{EF}, II, we introduced commuting difference analogues of the 
differential operators provided by action of the Sugawara tensor at 
critical level. It would be interesting to find commuting difference analogues
of the operators $T^{diff}(z)$ of the present work. 

4) In sect.\ \ref{sect:filt}, we introduce filtrations of the  conformal
blocks, provided by the maps $\cH_\uP$; it would be interesting  to
compute the corresponding $q$-dimensions. 

5) One could expect simple transformations from our version of the KZB
system to the versions of van Geemen-Previato (\cite{Previato}) and of 
Gaw\c{e}dzki-Tran (\cite{Gaw}), when the genus of $C$ equals $2$.


\subsection{}

The first author would like to express his thanks to G.\ Felder for many
conversations on conformal blocks and collaboration in 
\cite{EF,EF:sols}.  We also would like to thank V.\ Drinfeld, whose talk
at the 1994 Alushta conference (Ukraine) introduced us to the subject 
of this work. We also would like to thank  A.\ Lossev and E.\ Vasserot
for discussions related to this work. 

We also would like to express our gratitute to A.\ Alekseev, L.\ 
Faddeev and H.\ Grosse for an invitation at ESI (Vienna) in August 
1999, during which this work was started. The second author was partially
supported by grant RFFI-98-01-00327.

\section{Hecke parametrization of $\Bun_2(C)_{\cO(\sum_i P_i)}$} 
\label{sect:1}

Let $C$ be a smooth, complete, connected complex curve, let $g$ be its
genus. Let us set  $G = SL_2$, so $G(\CC) = SL_2(\CC)$ . For $\cL$ a
line bundle over $C$, let us denote  by $\Bun_2(C)_\cL$ the moduli space
of rank $2$ vector bundles over $C$, whose determinant bundle is
isomorphic to $\cL$.

\subsection{Theta-functions} \label{sect:theta}

Let us fix our conventions for theta-functions. We set, for
$\la$ in $\CC^g$, and $\tau = (\tau_{ab})_{a,b = 1,\ldots,g}$
a symmetric matrix with negative real part, 
$$
\Theta(\la |\tau) = \sum_{m\in\ZZ^g} \exp({1\over 2} m \tau m^t + 
m \la^t) . 
$$
We set $\tau(C)_{ab} = \int_{B_a}\omega_b$, and 
 $L(C) = (2i\pi \ZZ)^g + \tau\ZZ^g$. Then the degree zero part $J^0(C)$
of the Jacobian of $C$ is isomorphic to $\CC^g / L(C)$. The Abel-Jacobi map
sends a divisor $\sum_i n_i P_i$ of $C$ of total degree zero to 
the class of the vector 
$(\sum_i \int_{x_0}^{P_i} \omega_a)_{a = 1,\ldots, g}$ (which 
is independent of the choice of $x_0$ and of the integration contours). 

The vector of Riemann constants $\Delta$ is the element of the
degree $g-1$ part $J^{g-1}(C)$ of the Jacobian of $C$ such that the 
identity $\Theta(\sum_{i = 1}^{g-1} x_i - \Delta) = 0$. We also 
fix an odd theta-cahracteristic $\delta$ (the product of $\la\mapsto
\Theta(\la + \delta - \Delta)$ by an exponential factors is then 
an odd function of $\la$).  

Let us set 
\begin{equation} \label{r:P}
r^{(P)}(z) = d_P \ln\Theta( A(P) - A(z) + \delta - \Delta).
\end{equation}
Then we have 
$$
r^{(P)}(\gamma_{A_a}z) = r^{(P)}(z), \quad  r^{(P)}(\gamma_{B_a}z) 
= r^{(P)}(z) + \omega_a(P) . 
$$

\subsection{The Hecke parametrization map $\cH_{\uP}$}

Let $\cE$ be a vector bundle of rank $2$ over $C$. Let  $\uP = (P_i)_{i
= 1, \ldots, N}$ be a family of points on $C$. For  $P$ a point of $C$,
let us denote by $\cE_P$ the fiber of  $\cE$ at $P$. Let $\ul =
(\ell_i)_{i = 1, \ldots, N}$ be an element of  $\prod_{i = 1}^N
\PP(\cE_{P_i})$. The {\it Hecke modification} $\cH_{\uP,\ul}(\cE)$ of
$\cE$ along $(\uP,\ul) = (P_i,\ell_i)_{i = 1,\ldots,N} $ is defined as
the sheaf whose space of sections over an open subset $U$ of $C$  is
$\{\sigma:$ rational section of $\cE$ over $U$, regular   outside
$\{P_i, i = 1,\ldots,N\}$, with a simple pole at each $P_i$ and residue
in the line $\ell_i\}$. 

Define a map $\cH_{\uP} : (\CC P^1)^N \to \Bun_2(C)_{\cO(\sum_i P_i)}$
by 
$$
\cH_{\uP}(\ul) := \cH_{\uP,\ul}(\cO^2_C), 
$$
where $\ul = (\ell_1,\ldots,\ell_N)$ belongs to $(\CC P^1)^N$ and 
$\cO^2_C$ denotes the trivial bundle of rank $2$ over $C$.  Since the
group of automorphisms of $\cO^2_C$ is isomorphic to  $GL_2(\CC)$,
$\cH_{\uP}$ factors through $(\CC P^1)^N \to (\CC P^1)^N /GL_2(\CC) = 
(\CC P^1)^N /G(\CC)$.

\subsection{Properties of $\cH_\uP$ in the case $N = 3g$ (proof of 
Thm.\ \ref{thm:local:isom})}

\subsubsection{The preimage of stable bundles (proof of the first 
part of Thm.\ \ref{thm:local:isom})} \label{AJ}

In this section, we will characterize the preimage by the map $\cH_{\uP}$ of
the subset of $\Bun_2(C)_{\cO(\sum_i P_i)}$ consisting of stable bundles. 

Let us define $M(\uP,\ul)$ as the $3g \times 3g$ matrix with
matrix elements
$$
M(\uP,\ul)_{i,3j+\al} = \omega_{a}(P_{i}) \ell_{i}^{\al}.
$$
for $i = 1,\ldots,3g,a = 1,\ldots,g, \al = 1,2,3$.

\begin{lemma}
We have 
$$
\det M(\uP,\ul) = \Den(\uP,\ul) ,  
$$
where $\Den(\uP,\ul)$ is defined by (\ref{expr:Den}).  
\end{lemma}

Recall that the bundle $\cE$ is called stable iff $H^0(C,\End(\cE)) =
\CC$. The set of stable bundles has the structure of  a quasi-projective
variety (\cite{Seshadri}).  

\begin{prop}  \label{prop:stable}
Assume that all $\ell_i$ are  $\neq \infty$. 
Then if  $\Den(\uP,\ul) \neq 0$, $\cH_{\uP}(\ul)$
is stable.  
\end{prop}

{\em Proof.} Let us set $\cE = \cH_{\uP}(\ul)$. 
Let $\End(\cE)_0$ be the subbundle of $\End(\cE)$ formed of traceless
endomorphisms. We derive from the exact sequence of sheaves
$$
0 \to \End(\cE)_0 \to \End(\cE) \to \cO_C \to 0, 
$$
the exact sequence $0 \to H^0(C,\End(\cE)_0)   
\to H^0(C,\End(\cE)) \to \CC$. The stability  condition is therefore 
equivalent to $H^0(C,\End(\cE)_0)  = 0$. 

Let us denote by $\CC(C)$ the function field of $C$.
$H^0(C,\End(\cE)_0)$ is isomorphic to the subspace of  $\bar\G \otimes \CC(C)$
formed of the rational functions $\varphi : C\to \bar\G$, 

{\it i)} regular outside  $\{P_i\}$, 

{\it ii)} which have simple poles at each $P_i$ and residue at that point 
in $\CC( - e + \ell_i h + \ell_i^2 f)$, 

{\it iii)} such that if we write $\varphi = \varphi_e e  + \varphi_h h  +
\varphi_f f$, with $\varphi_e,\varphi_h,\varphi_f$ in $\CC(C)$,  $ -
\ell_i^2 \varphi_e  - 2 \ell_i \varphi_h + \varphi_f$ is regular at
$P_i$.  

Conditions {\it i)} and {\it ii)} mean 
that for some constants $(\la_i)_{i = 1,\ldots, 3g}$ and $C_e,C_h,C_f$, we have 
$$
\varphi = \sum_{i = 1}^{3g} \la_i ( - e + \ell_i h + \ell_i^2 f) r^{(P_i)}(z)
+ C_e e + C_h h + C_f f 
$$
and
\begin{equation} \label{system:1}
\sum_i \la_i \ell_i^\al\omega_a(P_i) = 0, \quad a = 1,\ldots, g, \quad 
\al = 0,1,2 ,   
\end{equation}
and condition {\it iii)} is equivalent to 
\begin{equation} \label{iii}
\sum_{j\neq i} \la_j \ell_{ij}^2 r^{(P_j)}(P_i) +  C_f - 2 \ell_i C_h - \ell_i^2 C_e = 0  
\on{\ for\  each\ } i = 1, \ldots, 3g.  
\end{equation}
Since $\Den(\uP,\ul)\neq 0$, (\ref{system:1}) implies that 
all $\la_i$ are zero. Moreover, $\Den(\uP,\ul)\neq 0$ also implies that  
card$\{\ell_i; i = 1,\ldots, 3g\}$ is $\geq 3$. Then 
(\ref{iii}) implies that $C_f = C_h = C_e = 0$. \hfill \qed\medskip  

\subsubsection{The fibers of $\cH_{\uP}$}

In this section, we prove: 

\begin{prop} \label{prop:2.2}
For any $\ul$ such that $\Den(\uP,\ul) \neq 0$, the set 
$\cH_\uP^{-1}(\cH_{\uP}(\ul))$ is finite. 
\end{prop}

{\em Proof.} Let $\ul$ and $\ul'$ belong to $(\CC P^1)^{3g}$. Assume that $\Den(\uP,\ul)
\neq 0$. The bundles $\cH_\uP(\ul)$ and 
$\cH_\uP(\ul')$ are isomorphic iff there exists a rational function $M(z)$ defined on 
$C$, with values in $SL_2(\CC)$, regular except for simple poles at the $P_i$, 
where the residue is proportional to $\pmatrix 1 \\ \ell'_i\endpmatrix 
\pmatrix \ell_i & -1 \endpmatrix$. $M(z)$ is then written 
$$
M(z) = C + \sum_i \al_i  \pmatrix 1 \\ \ell'_i\endpmatrix 
\pmatrix \ell_i & -1 \endpmatrix r^{(P_i)}(z) , 
$$
where $C$ is in $M_2(\CC)$ and  
where $\al_i$ are tangent vectors at $P_i$, satisfying 
\begin{equation} \label{unival}
\sum_i \omega_a(P_i) \al_i  \pmatrix 1 \\ \ell'_i\endpmatrix 
\pmatrix \ell_i & -1 \endpmatrix = 0, \quad a=  1,\ldots, g, 
\end{equation}
and
\begin{equation} \label{van:residue}
\pmatrix \ell'_i & -1 \endpmatrix C \pmatrix 1 \\ \ell_i \endpmatrix 
+ \sum_{j\neq i} r^{(P_j)}(P_i) \al_j (\ell'_i - \ell'_j)  (\ell_j - \ell_i) = 0 , 
\quad i=  1,\ldots, 3g.   
\end{equation}
Let us set $\beta_i = \al_i \ell'_i$.  Since $\Den(\uP,\ul) \neq 0$, 
(\ref{unival}) means that the vector
$(\al_i,\beta_i)_{i = 1,\ldots, 3g}$ belongs to a $2g$-dimensional vector space. 
On the other hand, (\ref{van:residue}) is equivalent to 
\begin{equation} \label{quadratic}
\pmatrix \beta_i & -\al_i \endpmatrix C \pmatrix 1 \\ \ell_i \endpmatrix 
+ \sum_{j\neq i} r^{(P_j)}(P_i) 
(\al_j \beta_i - \al_i \beta_j)  (\ell_j - \ell_i) = 0, 
\quad i=  1,\ldots, 3g , 
\end{equation}
which is a system of $3g$ quadratic conditions on this vector. 
If these conditions are satisfied, the condition that 
$\det M(z) = 1$ is then equivalent to $\det M(P_0) = 1$. 
The set of $(C,\al_i,\beta_i)$ satisfying these conditions is 
therefore  a subvariety of  
$M_2(\CC) \times \CC^{3g} \times \CC^{3g}$, and 
$\cH_\uP (\cH_\uP^{-1}(\ul))$ is the image of this variety by the 
morphism $(C,\al_i,\beta_i)_{i = 1,\ldots, 3g} \mapsto 
(\beta_i / \al_i)_{i = 1,\ldots, 3g}$. 

Let us fix now $\ul'$ in $\cH_\uP (\cH_\uP^{-1}(\ul))$. It
follows from Prop.\ \ref{prop:stable} that $\Den(\uP,\ul)\neq 0$. 
Let us compute the tangent space of $\cH_\uP (\cH_\uP^{-1}(\ul))$ at
$\ul'$. Let $\delta \ul'$ be an infinitesimal element of this 
vector space. To it is associated a matrix 
corresponding to the isomorphism from 
$\cH(\uP,\ul')$ to $\cH(\uP',\ul' + \delta\ul')$. We have 
$$
M(z) = \id_{\CC^2} + \sum_i 
\al_i \pmatrix 1 \\ \ell'_i + \delta\ell'_i\endpmatrix  
\pmatrix \ell'_i & -1
\endpmatrix  r^{(P_i)}(z) , 
$$
where $\al_i$ are infinitesimals. They should satisfy 
the equations
\begin{equation} \label{urgences}
\sum_i \al_i \omega_a(P_i) \pmatrix 1 \\ \ell'_i + \delta\ell'_i \endpmatrix 
\pmatrix \ell'_i & -1 \endpmatrix = 0  
\end{equation}
for each $a$, therefore 
$$
\sum_i \al_i \omega_a(P_i)  = \sum_i \al_i \ell'_i\omega_a(P_i)   
= \sum_i \al_i \ell^{\prime 2}_i\omega_a(P_i)   = 0. 
$$
As $\Den(\uP,\ul') \neq 0$, this implies that the $\al_i$ are all
zero. 
\hfill \qed \medskip

\begin{remark}
The subvariety of $M_2(\CC) \times \CC^{3g} \times \CC^{3g}$
defined by conditions (\ref{unival}) and (\ref{van:residue})   (without
the condition $\det M(P_0) = 1$) has components of dimension $>0$. This
makes it difficult to evaluate  the degree of $\cH_\uP$. 
\end{remark}

\subsubsection{$\cH_{\uP}$ is a finite etale morphism (proof of the 
second part of Thm.\ \ref{thm:local:isom})} \label{sect:etale}

In \cite{Seshadri}, $\Bun_2(C)_{\cO(\sum_i P_i)}^{stable}$  is endowed
with the structure of a quasi-projective  variety. On the other hand,
since the action of $G$ on  $\{\ul \in (\CC P^1)^{3g} | 
\Den(\uP,\ul)\neq 0\}$ is free,  the quotient  $\{\ul \in (\CC P^1)^{3g}
|  \Den(\uP,\ul)\neq 0\} / G$  also has the structure of a 
quasi-projective variety. One can check using  \cite{Seshadri} that the
restriction of $\cH_\uP$ to  $\{\ul \in (\CC P^1)^{3g} | 
\Den(\uP,\ul)\neq 0\} / G$  is a morphism of varieties. 

Moreover, the dimensions of  $\{\ul \in (\CC P^1)^{3g} | 
\Den(\uP,\ul)\neq 0\} / G$  and  $\Bun_2(C)_{\cO(\sum_i P_i)}^{stable}$
are equal. It follows from Prop.\ \ref{prop:2.2} that the kernel of the
map induced by $\cH_\uP$ on tangent spaces is zero. It follows  that
$\cH_\uP$ induces an isomorphism between the tangent spaces of $\{\ul
\in (\CC P^1)^{3g} |  \Den(\uP,\ul)\neq 0\} / G$  at $G\ul$  and  of
$\Bun_2(C)_{\cO(\sum_i P_i)}^{stable}$ at $\cH_\uP(\ul)$.  It follows
that the restriction $\cH_\uP$ to  $\{\ul \in (\CC P^1)^{3g} | 
\Den(\uP,\ul)\neq 0\} / G$  is etale. \hfill \qed\medskip

\begin{remark} The map $\cH_\uP$ is never surjective; for example, 
the $p$ large enough, the bundle $\cO(pP_0) \oplus \cO(\sum_i P_i 
- p P_0)$ is not in the image of $\cH_\uP$. But one may ask whether
$\Imm \cH_\uP$ contains all stable bundles when 
$N = \on{card}\{P_i\}$ is large enough. 
\end{remark}

\section{The Hitchin system in Hecke coordinates} \label{sect:hitchin}

\subsection{The map induced by $\cH_{\uP}$ on cotangent bundles} 

Let us assume that $N = 3g$. Let $\ul = (\ell_i)_{i = 1,\ldots,3g}$ be 
such that $\Den(\uP,\ul) \neq 0$. It follows from Thm.\
\ref{thm:local:isom} that  $\cH_{\uP}$ is an etale isomorphism 
at the neighborhood of  the class $[\ul]$ of $\ul$ in $(\CC P^1)^{3g} / G(\CC)$.

$\cH_{\uP}$ therefore induces an isomorphism of cotangent spaces 
$$
T^* \cH_{\uP} : T^*_{[\ul]} [(\CC P^1)^{3g} / G(\CC)]
\to T^*_{\cH_\uP(\ul)} \Bun_2(C)_{\cO(\sum_{i = 1}^{3g} P_i)} . 
$$

\subsubsection{} \label{identification:cotgt:CP1}
For $\ell$ in $\CC P^1$, $T^*_\ell\CC P^1$ is
isomorphic to $\CC$.  Let us denote by $\bar\G$ the Lie algebra
$\SL_2(\CC)$. Let $e,h,f$ be its Chevalley generators, and 
let $\langle , \rangle_{\bar\G}$ be the invariant form on 
$\bar\G$ such that $\langle e,f \rangle_{\bar\G} 
= {1\over 2} \langle h,h \rangle_{\bar\G} = 1$. We will identify 
$\bar\G$ with its dual using $\langle , \rangle_{\bar\G}$. 
$T^*\CC P^1$ is a symplectic manifold with a Hamiltonian  action of
$SL_2(\CC)$, and the corresponding moment map  $\mu : T^* \CC P^1
\to \bar\G^*$ is given by $(\ell,\la) \mapsto \la (-e + \ell h + \ell^2
f)$;  $\mu$ coincides with the Springer desingularization. 

$T^*_{[\ul]} [(\CC P^1)^{3g} / G(\CC)]$ is therefore isomorphic
to $(\mu^{3g})^{-1}(0) / G(\CC)$; the preimage of $\CC^{3g}$ by the projection 
$\mu^{-1}(0) \to (\CC P^1)^{3g}$ is 
$$
\{ (\ul,\ula) = (\ell_i,\la_i)_{i = 1, \ldots, 3g} | \sum_{i = 1}^{3g} \la_i = 
\sum_{i = 1}^{3g} \la_i \ell_i  = \sum_{i = 1}^{3g} \la_i \ell_i^2 = 0 \} , 
$$
and the action of $G(\CC)$ on this space is 
$$
\pmatrix a & b \\ c & d \endpmatrix \cdot (\ell_i,\la_i)_{i = 1,\ldots,3g} 
= ({{ a\ell_i + b}\over{ c\ell_i + d}} , {{\la_i}\over{ (c\ell_i + d)^2 }}
)_{i = 1,\ldots,3g} . 
$$

\subsubsection{} Let $\cE$ be any bundle of $\Bun_2(C)_{\cO(\sum_{i =
1}^{3g} P_i)}$.  $T^*_\cE \Bun_2(C)_{\cO(\sum_{i = 1}^{3g} P_i)}$ is
equal to $H^0(C,\End(\cE)_0  \otimes \Omega_C)$, where $\Omega_C$ is
the sheaf of differentials on $C$ and  $\End(\cE)_0$ is the sheaf of
traceless endomorphisms of $\cE$.

\begin{lemma} \label{lemma:cotgt:moduli}
Let $\cE = \cH_{\uP}(\ul)$. Let us define $\Omega(C)$ as the space
of rational differentials on $C$.  Then $H^0(C,\End(\cE)_0 \otimes
\Omega_C)$ is isomorphic with the subspace of $\Omega(C)
\otimes \SL_2$ consisting of the elements $A(z)$ such that

i) $A(z)$ is regular on $C\setminus \{P_i, i = 1,\ldots, 3g\}$ 
and has simple poles at the $P_i$, 

ii) for any $i = 1, \ldots, 3g$,  $\res_{z = P_i} A(z)$ belongs to $\CC(-e +
\ell_i h + \ell_i^2 f)$, 

iii) if we set $A(z) = A_e \otimes e + A_h \otimes h + A_f \otimes f$, with
$A_e,A_h,A_f$ in $\Omega(C)$, then $-\ell_i^2 A_e - 2 \ell_i A_h +
A_f$ is regular at $P_i$, for any $i = 1, \ldots, 3g$.  
\end{lemma}

Condition {\it iii)} means that if $\rho$ is the fundamental
representation of  $\SL_2$, the form $\det(\pmatrix 1\\ \ell_i
\endpmatrix , (\id\otimes \rho) (A(z)) \pmatrix 1\\ \ell_i \endpmatrix)$ is
regular at $P_i$.

\subsubsection{}

Let us fix a basis $(A_a,B_a)_{a = 1,\ldots,g}$ of $a$- and $b$-cycles
on $C$ and let $(\omega_a)_{a = 1,\ldots,g}$ be the associated Abelian
differentials on $C$.  

Let us set $\omega^{(P)}(z) = r^{(z)}(P)$ (see formula (\ref{r:P})). 

\begin{lemma} 
In terms of the identifications of \ref{identification:cotgt:CP1} and
Lemma \ref{lemma:cotgt:moduli}, the map $T^*_{\ul} \cH_{\uP} $ is
given by
\begin{equation} \label{formula:A}
T^*_{\ul} \cH_{\uP}
(\la_1,\ldots,\la_{3g}) = A(\ul,\ula|z) , 
\end{equation}
where $A(\ul,\ula|z)$ is defined by (\ref{def:A})
\end{lemma}

\subsection{Expression of the Hitchin Hamiltonians}

\subsubsection{}

The Hitchin fibration is defined as the map 
$$
T^* \Bun_2(C)_{\cO(\sum_i P_i)} \to H^0(C,(\Omega_C)^{\otimes 2}) , 
$$
$(\cE,\xi)\mapsto (\id \otimes \tr)(\xi^2)$, where $\xi$ belongs to $H^0(C,\Omega_C
\otimes \End(\cE)_0)$. The Hitchin Hamiltonians are defined as follows. 
Let $\omega^{(2)}_1,\ldots,\omega^{(2)}_{3g-3}$ be a basis 
of the space of quadratic differentials on $C$ and let us set 
$$
(\id\otimes \tr)(\xi^2) = \sum_{\al = 1}^{3g-3} \Hitch_\al(\cE,\xi) \omega^{(2)}_\al ; 
$$
the $(\cE,\xi)\mapsto \Hitch_\al(\cE,\xi)$ form a Poisson-commutative
family of functions on  $T^* \Bun_2(C)_{\cO(\sum_i P_i)}$, the Hitchin
Hamiltonians (see \cite{Hitchin}).

\subsubsection{}

The right side of (\ref{formula:A}) is equal to $A(\ul,\ula|z)$ defined in 
(\ref{def:A}). 
Prop.\ \ref{prop:hitchin} follows from 
$$
H(\ul,\ula |z) = \tr \rho[A(\ul,\ula|z)]^2 ,  
$$
where $\rho$ is the fundamental representation of $\SL_2$. 

Another expression of $H(\ul,\ula |z)$
may be obtained as follows. Let us define $A^{reg}$ as the second sum of 
(\ref{formula:A}), and $\nu_{ki}(z)$ as $\det(a_{ki})$, where $a_{ki}$
is the matrix with coefficients
$(a_{ki})_{j,3 a - \al} = \ell_j^{\al} \omega_a(P_j)$ if $j\neq k$, and 
$(a_{ki})_{k,3a -\al} = \ell_i^\al \omega_a(z)$ if 
$a = 1,\ldots,g, \al = 1,2,3$. Then 
\begin{align*}
H(\ul,\ula |z) = \tr \rho(A^{reg})^2 & + 
{1\over {\Den(\uP,\ul)}} \sum_{i,j,k; k\neq i}
\la_i\la_j \ell_{ki}^2 \nu_{kj}(z) \omega^{(P_i)}(P_k) \omega^{(P_j)}(z) 
\\ & 
- {1\over 2} \sum_{i,j; i\neq j} \la_i \la_j \ell_{ij}^2  
\omega^{(P_i)}(z) \omega^{(P_j)}(z) . 
\end{align*}

\section{Parametrization of conformal blocks (the case $N = 3g$)}
\label{sect:cb}

\subsection{``Twisting'' group elements}

\subsubsection{Element $g_{\ul,\uP}$ at level zero} \label{sect:def:b}

Let $P_0$ be a point of $C$, distinct from the $P_1,\ldots, P_{3g}$.  
Let $J(C)$ be the Jacobian of $C$; $J(C)$ is the direct sum $\oplus_{i\in\ZZ} J^i(C)$
of its graded components.  Let us denote by $A : C \to J^1(C)$ the Abel-Jacobi map. 
We assume that 
\begin{equation} \label{ass:jac}
A(P_1) +\ldots + A(P_{3g}) = 3g A(P_0) . 
\end{equation} 

(The results of this section can easily be transposed 
when $3g P_0$ is replaced by an arbitrary effective divisor 
of degree $3g$.)  

Let us set, for $z$ in $C$, 
$$
b_\uP(z) = {{
\prod_{i = 1}^{3g} \Theta( A(z) - A(P_i) + \delta - \Delta) 
}\over{\Theta(A(z) - A(P_0) + \delta - \Delta)^{3g} }} . 
$$ 
Assumption (\ref{ass:jac}) implies that $b_\uP$ is a rational 
rational function on $C$, with divisor $\sum_{i = 1}^{3g} P_i - 3g P_0$
(so $b_\uP$ has zeroes at the $P_i$ and a pole at $P_0$). Let us fix 
functions
$d_{i,\uP}$, regular on $C\setminus \{P_0\}$, such that $d_{i,\uP}(P_j) = \delta_{ij}$. 
We then set, for $\ell_i$ in $\CC P^1\setminus 
\{0\}$, and $z$ in $C \setminus \{P_0,P_i, i = 1,\ldots, 3g\}$,   
$$
g_{\ul,\uP}(z) = \pmatrix b_\uP^{-1} & 0 \\ 0 & 1\endpmatrix \pmatrix 1 & -
\sum_i \ell_i^{-1} d_{i,\uP}\\ 0 & 1\endpmatrix(z) . 
$$

Recall that for $A$ a ring, $GL_2(A)$ is defined as the group of
matrices  $\pmatrix a & b \\ c & d \endpmatrix$, with $a,b,c,d$  in $A$,
 and $ad - bc$ in $A^\times$ ($A^\times$ denotes the set of 
invertible elements of $A$).  For each point $P$ of $C$, we fix a local
coordinate $z_P$ at $P$. Recall that $GL_2(\CC(C))$ may be viewed as a
subgroup of the adeles  group $GL_2(\AAA_C) = \prod'_{P\in C}
GL_2(\CC((z_P)))$.

\begin{lemma}  \label{lemma:pties:gPl} 
For each family $\ul = (\ell_i)_{i = 1,\ldots, 3g}$ of $(\CC P^1 \setminus 
\{0\})^{3g}$, the function $g_{\ul,\uP}$ satisfies the equivalent conditions

1) $z\mapsto g_{\ul,\uP}(z)$ is a rational function  on $C$ with 
values in $G$, regular on $C \setminus \{P_0;P_i,i = 1,\ldots, 3g\}$, 
with  simple poles at each $P_i$, $i = 1,\ldots, 3g$, the expansion
$$
g_{\ul,\uP}(z) = {1\over{z_{P_i}}} \pmatrix a_i \\ c_i \endpmatrix 
\pmatrix \ell_i & -1 \endpmatrix + O(1) 
$$ 
at $P_i$, and $\det(g_{\ul,\uP})$ has the expansion  $
\det(g_{\ul,\uP}) = c_i z_{P_i}^{-1} + O(1)$ at the same point,  where
 $a_i,b_i$ and $c_i$ are complex numbers, 
such that $(a_i,b_i)\neq (0,0)$  and $c_i\neq 0$; 

2) $g_{\ul,\uP}$ belongs to the intersection of $GL_2(\CC(C))$ and 
the subset 
$$ 
GL_2(\CC((z_{P_0})))
\times \prod_{i = 1}^{3g}  GL_2(\CC[[z_{P_i}]])
\on{diag}(1,z_{P_i}^{-1})  \pmatrix 1 & 0 \\ - \ell_i & 1\endpmatrix 
\times \prod_{Q\notin \{P_0,P_i \} } GL_2(\CC[[z_{Q}]])$$
of $GL_2(\AAA_C)$.  
\end{lemma}

The properties of $g_{\ul,\uP}(z)$ with respect to variation of $\ul$ 
can be described as follows. Let $\cK$ be the local field of
$C$ at $P_0$, so $\cK = \CC((z_{P_0}))$. Let us set $R = H^0(C \setminus
\{P_0\},\cO_C)$, which we view as a subring of $\cK$, and let us set
$\G^{out} = \bar\G\otimes R$. It is a Lie subalgebra of 
$L\bar\G = \bar\G\otimes \cK$. Let us denote by $\wt\G^{out}_{\ul}$ 
the subspace of $L\bar\G$ formed of the Laurent expansions at $P_0$
of the rational maps $\phi: C \to \bar\G$, regular outside $\{P_0,P_1,\ldots,P_{3g}\}$, 
with expansion at each $P_i,i\neq 0$
$$
\phi(z) = {{\la_i}\over{z_{P_i}}} (- e + \ell_i h + \ell_i^2 f) + O(1) , 
$$
where the $\la_i$ are complex numbers. 

Define for $i = 1,\ldots, 3g$, $\mu_i$ as the linear map from 
$\wt\G^{out}_{(\ell_i)}$ to $\CC$ given by 
$$
\mu_i(\phi) = ( - \phi_f  + 2 \ell_i\phi_h  + \ell_i^2\phi_e)(P_i) , 
$$
for $\phi = \phi_e e + \phi_h h + \phi_f f$, $\phi_e,\phi_h,\phi_f$ in 
$\CC(C)$. Let us set 
$\G^{out}_{\ul} := \Ad(g_{\ul,\uP})^{-1}(\G^{out})$
($\Ad$ denotes the adjoint action of a group on its Lie algebra). 
Then we have 
\begin{equation} \label{exact}
\G^{out}_{\ul}
= \Ker \left(\oplus_{i=1}^{3g} \mu_i: \wt \G^{out}_{\ul} \to \CC^{3g} \right) . 
\end{equation}
Moreover, we have 
\begin{lemma} \label{diff:eqs}
For any $\phi$ in $\wt\G^{out}_{\ul}$, we have 
$$
\Ad(g_{\ul,\uP})(\phi) \in \sum_{i = 1}^{3g} \mu_i(\phi) (\pa_{\ell_i} 
g_{\ul,\uP}) g_{\ul,\uP}^{-1} + \G^{out}
$$
\end{lemma}

{\em Proof.} One first shows this statement for $\phi$ of the form 
$e\otimes r$, $r\in R$, using the explicit formula for $g_{\ul,\uP}$. 
In the case of a general $\phi$, one then replaces $\phi$ by $\phi +
e\otimes \{\sum_i \ell_i^{-2}\mu_i(\phi) d_{i,\uP} \}$ and uses 
(\ref{exact}). \hfill \qed\medskip

\subsubsection{Kac-Moody algebras, extended loop groups and 
semidirect products} \label{sect:repr} \label{semidirect}

Let $\G$ be the Kac-Moody  Lie algebra $\G  = L\G \oplus \CC K$. It is
endowed with the Lie bracket $[(x\otimes a, bK), (x'\otimes a',b'K)]  =
( [x,x']\otimes aa', \langle x,x' \rangle_{\bar\G} \res_{P_0}(da\cdot
a') K)$. For $x$ in $\bar\G$ and $\phi$ in $\cK$, we will denote
$(x\otimes \phi,0)$ by $x[\phi]$.  
  
Recall that for any ring $A$, $G(A)  = SL_2(A)$ is the group of matrices 
$\pmatrix a & b \\ c & d \endpmatrix$, where $a,b,c,d$ belongs to $A$ and such that 
$ad - bc = 1$ . 

Let $k$ be an integer $\geq 0$, and let us denote by $\wh{G(\cK)}$   the
level $k$ central extension of $G(\cK)$  by $\CC^{\times}$  (see
\cite{Garland}).  We denote by $Z$ the center of $\wh{G(\cK)}$. $Z$
is isomorphic with $\CC^\times$; for $c$ in $\CC^\times$, we denote by $z(c)$
the corresponding element of $Z$. 

Let $k$ be an integer $\geq 0$, and let us denote by $\wh{G(\cK)}$   the
level $k$ central extension of $G'(\cK)$  by $\CC^{\times}$  (see
e.g.\ \cite{Garland,Kumar}).  We denote by $Z$ the center of $\wh{G(\cK)}$. $Z$
is isomorphic with $\CC^\times$; for $c$ in $\CC^\times$, we denote by $z(c)$
the corresponding element of $Z$. 

We will use the following properties of the central extension
$\wh{G(\cK)}$: 

a) if $(\rho_\VV,\VV)$ is any integrable $\G$-module of level  $k$,
there is a unique lift (denoted $\rho^{group}_\VV$) of the action of
$\G$ to a representation of $\wh{G(\cK)}$ on $\VV$;  this lift satisfies
$\rho_\VV^{group}(z(c)) = c \id_\VV$, for $c$ in $\CC^\times$, and the 
compatibility rule $\rho_\VV(\Ad(\pi(g))(x)) = 
\rho_\VV^{group}(g)\rho_\VV(x)\rho_\VV^{group}(g)^{-1}$, for any $g$
in $\wh{G(\cK)}$ and $x$ in $\G$, where $\pi$ is the canonical map from 
$\wh{G(\cK)}$ to $G(\cK)$ and denotes the adjoint action of $G(\cK)$ on $\G$;   

b) let $N$ be the positive unipotent subgroup of $G$, then there is a
unique lift $\iota$ of $N(\cK)$ to $\wh{G(\cK)}$, such that for each formal series 
$\rho$ of $\cK$, we have 
$$
\rho_\VV^{group}(\iota (\pmatrix 1 & \rho \\ 0 & 1 \endpmatrix) )
= \rho_\VV (\exp( e[\rho] ) ) .  
$$

The adjoint actions of $\wh{G(\cK)}$ on $\G$ and on itself  
factor through actions of $G(\cK)$ on $\G$ and on $\wh{G(\cK)}$.  
One can show that 
these actions extend to actions of $GL_2(\cK)$ on these sets, which we
denote by $\Ad$. Let us denote by $\wt T_\G$  and $\wt T_{\wh{G(\cK)}}$ 
the automorphisms of $\G$ and
of $\wh{G(\cK)}$ equal to $\Ad(\pmatrix z_{P_0} & 0 \\ 0 & 1 \endpmatrix )$. 
Then $\wt T_\G$ is given by 
$$
\wt T_\G [(e\otimes\phi,aK)] = (e\otimes z_{P_0}\phi, aK),  
\wt T_\G [(f\otimes\phi,aK)] = (f\otimes z^{-1}_{P_0}\phi, aK), 
$$
$$ 
\wt T_\G [(h\otimes\phi,aK)] = (h\otimes \phi, \{a + \res_{P_0}({{dz}\over z}\phi) \}K) , 
$$  
and $\wt T_{\wh{G(\cK)}}$ is the unique group automorphism of $\wh{G(\cK)}$
extending the automorphism $\Ad(\pmatrix z_{P_0} & 0 \\ 0 & 1 \endpmatrix)$ 
of $G(\cK)$.  Moreover, $\wt T_{\G}^2$ and $\wt T_{\wh{G(\cK)}}^2$ are both inner 
automorphisms, equal to $\Ad \pmatrix z_{P_0} & 0 \\ 0 & z_{P_0}^{-1}\endpmatrix$. 

Define now $\wh{G(\cK)}\rtimes\ZZ$ as the semidirect product of $\wh{G(\cK)}$ with the action of $\ZZ$ provided
by $\wt T_{\wh{G(\cK)}}$.  Precisely, $\wh{G(\cK)}\rtimes\ZZ$ is the set 
$\wh{G(\cK)}\times\ZZ$, endowed with the product 
$$
(g,n)(g',n') = (g \wt T_{\wh{G(\cK)}}^n(g'), n+n') . 
$$
We will write $T_{\wh{G(\cK)}}$ for the element $(e_{\wh{G(\cK)}},1)$, where 
$e_{\wh{G(\cK)}}$ is the neutral element of $\wh{G(\cK)}$, so  $(g,n)$
will be equal to $gT_{\wh{G(\cK)}}^n$.   The adjoint action of
$\wh{G(\cK)}$ on $\G$ extends to an action of  $\wh{G(\cK)}\rtimes\ZZ$
on $\G$ by Lie algebra automorphisms,  which we also denote by $\Ad$; it
is defined by $\Ad( gT_{\wh{G(\cK)}}^n)(x) =  \Ad(g)[ \wt
T_{\G}^n(x) ]$. 

Let us now construct representations of the group $\wh{G(\cK)} \rtimes\ZZ$. 
Let $(\rho_\VV,\VV)$ be an irreducible integrable representation of $\G$ of 
level $k$. Two possibilities occur: 

a) if $\rho_\VV\circ T_\G$ is equivalent to $\rho_\VV$, let 
$T_\VV$ be a linear automorphism of $\WW$ such that $\rho_\VV \circ T = 
\Ad(T_\VV) \circ \rho_\VV$. 
We set $\WW = \VV$ and let $\rho_\WW$ be the map from $\wh{G(\cK)} \rtimes\ZZ$
to $\Aut_\CC(\WW)$, defined by $\rho_\WW^{group}(g T_{\wh{G(\cK)}}^n)
= \rho_\VV^{group}(g) \circ T_{\VV}^n$; 

b) if $\rho_\VV\circ T_\G$ is not equivalent to $\rho_\VV$, we set 
$\WW = \VV \oplus \VV$, and we define $\rho_\WW$ as the map from $\wh{G(\cK)} \rtimes\ZZ$
to $\Aut_\CC(\WW)$, defined by 
$$
\rho_\WW^{group}(g T_{\wh{G(\cK)}}^n) :=
\pmatrix \rho_\VV^{group}(g) & 0 \\ 0 & (\rho_\VV^{group}\circ \wt T_{\wh{G(\cK)}}^{-1})(g)
\endpmatrix   
\pmatrix 0 & \rho_\VV^{group}(\diag(z_{P_0},z_{P_0}^{-1}))
\\ \id_\VV & 0 \endpmatrix^n .   
$$

In both cases, $(\WW,\rho_\WW^{group})$ is an irreducible representation of 
$\wh{G(\cK)} \rtimes\ZZ$. 

\subsubsection{Definition of $\wt g_{\ul,\uP}$}

Identify the function $z\mapsto b_\uP(z)$ (sect.\ \ref{sect:def:b}) 
with its image in $\cK^\times$.  We have $b_\uP = z_{P_0}^{-3g}
b_0$, where  $b_0$  belongs to $\CC[[z_{P_0}]]^\times$. $b_0$ has two 
opposite square roots, which belong to $\CC[[z_{P_0}]]^\times$. 
Let us fix one of them, which we denote by $b_0^{1/2}$. 
Let us fix a lift $\wt t$ in $\wh{G(\cK)}$ of $\pmatrix
(b_0^{1/2})^{-1} & 0 \\ 0 & b_0^{1/2} \endpmatrix$,  and let us set 
$\wt b_\uP = T_{\wh{G(\cK)}}^{3g} \wt t$.   Let us finally set, for
$\ul = (\ell_1,\ldots,\ell_{3g})$ in $(\CC^\times)^{3g}$,  
$$
\wt g_{\ul,\uP}(z) := \wt b_\uP \cdot \iota (\pmatrix  1  & 
- \sum_i \ell_i^{-1} d_{i,\uP} \\ 0 &  1 \endpmatrix ) 
\cdot z ((\prod_{i = 1}^{3g} \ell_i)^{k}) , 
$$  
where the rational function $d_{i,\uP}$ is identified with its image in $\cK$.  

In the previous section, we studied subspaces $\G^{out},\G^{out}_{\ul}$
and $\wt\G^{out}_{\ul}$ of $\bar\G\otimes \cK$. We will embed them 
in $\G$ using the following maps: 

a) $\G^{out}$ is mapped to $\G$ by the map $x\mapsto (x,0)$; 

b) $\wt\G^{out}_{\ul}$ is mapped to $\G$ by the map 
$$
i : \phi = \phi_e e + \phi_h h + \phi_f f \mapsto (\phi, 
\sum_i (\phi_h + \ell_i \phi_e)(P_i) K)  
$$

c) $\G^{out}_{\ul}$ is mapped to $\G$ by the restriction of the previous map. 

We will denote by $\G^{out,ext}, \wt\G^{out,ext}_{\ul}$ and 
$\G^{out,ext}_{\ul}$ the images of these maps. 

\begin{lemma} \label{ille}
$\G^{out,ext}$ and $\G^{out,ext}_{\ul}$ are Lie subalgebras of $\G$, 
isomorphic to $\G^{out}$ and $\G^{out}_{\ul}$ respectively.  
Moreover, we have 
$$
\G^{out,ext}_{\ul} = \Ad(\wt g_{\ul,\uP})^{-1} (\G^{out,ext}) . 
$$
\end{lemma}

The analogue of Lemma \ref{diff:eqs} is 
\begin{lemma} \label{juan}
For any $\phi$ in $\wt\G^{out}_{\ul}$, we have 
$$
\Ad(\wt g_{\ul,\uP}) (i(\phi)) \in \sum_{i = 1}^{3g} \mu_i(\phi)
(\pa_{\ell_i} \wt g_{\ul,\uP}) \wt g_{\ul,\uP}^{-1} 
+ \G^{out,ext} . 
$$
\end{lemma}

{\em Proof.} One again checks this when $\phi = e\otimes\phi$ for $\phi$
in $R$, and uses then Lemma \ref{ille}. 
\hfill \qed\medskip 

\subsubsection{Properties of $\wt g_{\ul,\uP}$}

Let us set $G^{out} = G(R)$. There is unique lift of $G^{out}$
to a subgroup of $\wh{G(\cK)}$, which we denote by $G^{out,ext}$,
corresponding to the Lie algebra inclusion $\G^{out,ext}\subset \G$. 
This lift is such that for any integrable representation 
$\VV$, and any $\G^{out,ext}$-invariant linear form $\psi$
on $\VV$, $\psi$ is also $G^{out,ext}$-invariant.

\begin{lemma} \label{depp}
The rational functions $h: (\CC P^1)^{3g}\to \wh{G(\cK)}$,  
$\ul = (\ell_i)_{i = 1,\ldots,3g}\mapsto h_{\ul}$ satisfying the
conclusion  of Lemma \ref{juan} for generic $\ul$ are exactly the
functions  $h_{\ul} = z_0\cdot \gamma_{\ul}\cdot \wt
g_{\ul,\uP}$, where  $z_0$ is a (constant in $\ul$) element of the
center $Z$ of $\wh{G(\cK)}$,  and $\gamma_{\ul}$ is a rational
function from  $(\CC P^1)^{3g}$ to $G^{out,ext}$.  
\end{lemma}

\begin{corollary} \label{invariance}
Let us set, for $g = \pmatrix a & b \\ c & d \endpmatrix\in G(\CC)$
and $\ell$ in $\CC P^1$, $g\cdot \ell = {{a\ell + b}\over{c\ell + d}}$;
$(g,\ell)\mapsto g\cdot \ell$ is the homographic action of $G(\CC)$ on $\CC
P^1$. Let us denote by $(g,\ul) \mapsto g\cdot \ul$ the diagonal action of 
$G(\CC)$ on $(\CC P^1)^N$. 
For $g$ in $G(\CC)$, let us denote by $\wt g$ the image of $g$
by the lift map from $G^{out}$ to $G^{out,ext}$. Then for any $g$ in $G(\CC)$,
there exists a rational function
$\ul \mapsto \gamma(g,\ul)$ with values in $G^{out,ext}$, such
that $\wt g_{g \cdot \ul,\uP} = \gamma(g,\ul) \cdot \wt
g_{\ul,\uP} \cdot \wt g \cdot  z(\prod_{i = 1}^{3g} (c\ell_i + d)^{-k})$. 
\end{corollary}

{\em Proof.} Let us set $h_{g,\ul} :=  \wt g_{g \cdot \ul,\uP}
\cdot  \wt g^{-1} \cdot  z(\prod_{i = 1}^{3g} (c\ell_i + d)^{k})$. 
Then one checks that $h_{g,\ul}$ satisfies the conclusion 
of Lemma \ref{juan}. Then  Lemma \ref{depp} implies that  there exists an element  $z_g$
of $Z$ and a rational function $\ul\mapsto \gamma(g,\ul)$
with values in $G^{out,ext}$, such that $h_{g,\ul} =
z_g \cdot \gamma(g,\ul) \cdot \wt g_{\ul,\uP}$.  $g\mapsto
z_g$ is therefore a group homomorphism from $G(\CC)$ to $Z =
\CC^\times$, and is therefore trivial. \hfill \qed\medskip 

\begin{remark} It is easy to check directly Lemma \ref{invariance} when 
$g$ is lower-triangular. 
\end{remark}

\subsection{Conformal blocks and associated functions}
\label{sect:def:weyl} \label{rep:semidirect}

Let $k$ be  an integer $\geq 0$ and let $\WW$ be one of the level $k$
representations of $\wh{G(\cK)} \rtimes\ZZ$ constructed in sect.\
\ref{sect:repr}.  For any integer $a\geq 0$, let $(\rho_a,V_a)$ be the
irreducible $(a+1)$-dimensional representation of $\bar\G$. The
irreducible integrable representations of level $k$ are the $\VV_{a,k}$,
where $0\leq a\leq k$. $\VV_{a,k}$ may be  constructed as the quotient
of the Weyl module $U\G\otimes_{U\G_+} V_{a}$,  where $\G_+ = \bar\G
\otimes \CC[[t]] \oplus\CC K$ and $V_a$ is endowed with the 
$\G_+$-module structure defined by $\rho_{a,+}: \G_+ \to \End(V_a)$ such
that $\rho_{a,+}(x\otimes f, aK) = \rho_a(x) f(0) + ka\id_{V_a}$.  
  
So 

a) either $k$ is even and $\WW = \VV_{k/2,k}$; 

b) or we have, for $a\neq b$ and $a+b = k$, $\WW = \VV_{a,k} \oplus 
\VV_{b,k}$. 

We denote by $W$ the subspace of $\WW$ equal to  the images of the
subspaces of the Weyl modules $1\otimes V_{k/2}$ in case a), and
$(1\otimes V_a) \oplus (1\otimes V_b)$ in case b); so $W$ is isomorphic to
$V_{k/2}$ in case a), and to $V_a \oplus V_b$ in case b). In each case,
$W$ is equal to the subspace of $\WW$ formed of the elements annihilated by
$\G^{in}$.  

Let $\psi$ be a $\G^{out}$-invariant linear form on $\WW$. We associate to 
it the function 
$$
f_\psi(\ul,\uP|v) := \langle \psi , 
\rho_\WW^{group}(\wt g_{\ul,\uP}) v\rangle , 
$$
where $v$ is in $W$. For any $\uP$,  $f_\psi(\ul,\uP|v)$ belongs
therefore to $\otimes_{i=1}^k \{\ell_i^k\CC[[\ell_i^{-1}]] \} \otimes
W^*$. 

Let us denote by $\CC[\ell]_{\leq k}$ the space of polynomials in $\ell$ 
of degree $\leq k$. This space is endowed with the right $G(\CC)$-module structure
defined by $(g\cdot f)(\ell) := f(g\cdot \ell) (cz + d)^{k}$.

{\em Proof of Prop.\ \ref{prop:fun}.} It follows from its definition that the map $v\mapsto
f_\psi(\ul|v)$ belongs to $\Hom_\CC(W,\otimes_{i = 1}^{3g} \ell_i^{k}
\CC[[\ell_i^{-1}]])$. Let may assume that $v$ is a weight vector,  so
$h[1]v = \la v$, where $\la$ is an integer in $\{-k,\ldots,k\}$. 
Moreover, 
$$
f_\psi(\ul,\uP|v) = \prod_{i=1}^g \ell_i^k
\sum_{\al\geq 0} {1\over{\al!}} \langle \psi, \rho_\WW^{group}(\wt b_\uP)
\rho_\WW(e[\sum_i d_{i,\uP}\ell_i^{-1}])^\al v \rangle.
$$
Since $\rho_\WW(h[1]) \rho_\VV^{group}(\wt b_\uP) = 
\rho_\WW^{group}(\wt b_\uP) [\rho_\WW(h[1]) - 3 gk]$, all contributions 
to this series where $\al \neq {{3gk + \la}\over 2}$ are zero. It follows that 
$$
f_\psi(\ul,\uP|v) = {1\over{\al!}} \langle \psi, \rho_\WW^{group}(\wt b_\uP)
\rho_\WW(e[\sum_i d_{i,\uP}\ell_i^{-1}])^\al v \rangle, \quad 
\on{with} \quad  \al = {{3gk + \la}\over 2} . 
$$
Therefore, for any $v$ in $W$,  $f_\psi(\ul,\uP|v)$ is a Laurent polynomial
in the $\ell_i$, of degree in each $\ell_i$ smaller than $k$. 

Finally, Cor.\ \ref{invariance}, together with the fact that
$\rho_\WW^{group}(\wt g) v = \rho_W(g)v$ for $g\in G(\CC)$, imply
that   $f_\psi( g\cdot \ul,\uP|v) =  \prod_{i=1}^{3g}
(c\ell_i + d)^k f_\psi(\ul,\uP|\rho_W(g)v)$, for  $g = \pmatrix a & b \\ c &
d\endpmatrix \in G(\CC)$. This formula both shows that $v\mapsto f_\psi(\ul,\uP|v)$  
is polynomial in the $\ell_i$, and that it is a morphism of $\bar\G$-modules. 
\hfill \qed\medskip 

We also get 

\begin{prop} 
If $\WW$ is a representation described of type b) above,  and 
if $\psi$ is zero when restricted to $\VV_{b,k}$, $f_\psi(\ul,\uP|v)$ is zero when 
$v$ belongs to $V_a$ if the genus $g$ of $C$ is odd, and when 
$v$ belongs to $V_b$ if $g$ is even. 
\end{prop}

\begin{remark}
The fact that $\WW$ is integrable should imply some 
additional functional properties of $f_\psi(\ul,\uP|v)$. 
\end{remark}

\begin{remark}
It follows from Lemma \ref{depp} that if we replace, in the definition
of $f_\psi(\ul,\uP|v)$, $\wt g_{\ul,\uP}$ by any function $h_{\ul,\uP}$  
satisfing the conclusions of Lemma \ref{juan}, the function
$f_\psi(\ul,\uP|v)$ only gets multiplied by a nonzero function of $\uP$. 
Moreover, since for any irreducible integrable representation
$\VV$, and any lift of $\pmatrix - 1 & 0 \\ 0 & -1 \endpmatrix \in
G(\cK)$, its image by $\rho_\VV^{group}$ is scalar,   a change in the
choice of the square root $b_0^{1/2}$ of $b_0$ only has the effect of 
multiplying $f_\psi(\ul,\uP|v)$ by a nonzero function of $\uP$. 
\end{remark}

\section{Action of the Sugawara tensor}

\subsection{Expressions of the Sugawara tensor $T(z)$}

For $x$ in $\bar\G$ and $\phi$ in $\cK$, let us denote by 
$x[\phi]$ the element $(x\otimes\phi,0)$ of $\G$. 

Let $\Cas$ be the element of  $\bar\G\otimes\bar\G$ corresponding to
$\langle, \rangle_{\bar\G}$. We have  $\Cas = {1\over 2}h\otimes h +
e\otimes f + f\otimes e$.  

Let us set for $x = e,h,f$, $x(z) = \sum_{i\in\ZZ} x[z_{P_0}^i]
z^{-i-1} dz$. 
We define the Green function $G(z,w)dz$ as 
$$
G(z,w)dz = d_z \ln\Theta(A(z) - A(w) + \delta - \Delta) 
- d_z \ln\Theta( A(z) - A(P_0) + \delta - \Delta). 
$$
We have then the expansion at the diagonal
$$
G(z,w)dz = {{dz}\over{z-w}} + O(1). 
$$ 
Expand $G(z,w)dz$ as $\sum_{i>0} \al_i(z) w^i$ for $w\ll z$ and 
as $\sum_{i>0} z^{i-1} dz \beta_i(w)$ for $z\ll w$ (both $z$ and $w$ being 
at the formal neighborhood of $P_0$). 
Let us then set 
$$
x^{in}(z) = \sum_{i>0} x[w^i] \al_i(z) , \quad  
x^{out}(z) = \sum_{i>0} x[\beta_i] z^{i-1} dz 
$$ 
($x^{in}(z)$ and $x^{out}(z)$ therefore belong to $\G\hat\otimes \Omega_z$)
and set 
\begin{equation} \label{def:sug}
T(z) = {1\over{2\kappa}} \sum_{\al} x_\al^{out}(z) y_\al(z) + y_\al(z)
x_\al^{in}(z) ,   
\end{equation}
where $\kappa = k + 2$, and $\sum_\al x_\al \otimes y_\al = \Cas$. 
The Laurent coefficients of $T(z)$  belong to a completion of the quotient
$U\G / (K - k\cdot 1)$ of the enveloping algebra of $\G$. 

The following Lemma shows that $T(z)$ is a normalization of the Sugawara
tensor.   
\begin{lemma} 
In any highest weight $\G$-module, the matrix elements 
of $T(z)$ coincide with those of 
$$
{1\over {2\kappa}}\limm_{z'\to z} \left( \sum_\al x_\al(z) y_\al(z')
- 3 d_{z'}[G(z,z')dz]  \right) .  
$$ 
\end{lemma}

We may also write (\ref{def:sug}) as follows. Let us set 
\begin{equation} \label{def:ell}
\ell^{in,std}(z) = \sum_\al x_\al \otimes y^{in}_\al(z), \ 
\ell^{out,std}(z)  = \sum_\al x_\al \otimes y^{out}_\al(z) ,  \
\ell(z) = \ell^{in,std}(z) + \ell^{out,std}(z). 
\end{equation}
Then 
$$
T(z) = {1\over{2\kappa}}
[(\tr\circ \rho)\otimes \id]\{ \ell^{out,std}(z) \ell(z) + \ell(z)
\ell^{in,std}(z) \} .  
$$
Let $\la(z) = \sum_\al x_\al \otimes \la_\al(z)$ be an element of 
$(L\bar\G\otimes \Omega) \otimes \bar\G$ and let us set now 
$$
\ell^{in}_\ul(z) = \ell^{in,std}(z) + \la(z), \quad
\ell^{out}_\ul(z) = \ell^{out,std}(z) - \la(z) . 
$$
Then 
\begin{align} \label{new:sug}
T(z) & = \nonumber 
{1\over{2\kappa}} [(\tr\circ \rho) \otimes \id ]
\{ \ell^{out}_\ul(z) \ell(z) + \ell(z) \ell^{in}_\ul(z) \} - [\la(z), \ell(z)] \} 
\\ & 
= {1\over{2\kappa}} \left( [\id \otimes (\tr\circ \rho)]\{ \ell^{out}_\ul(z) \ell(z) 
+ \ell(z) \ell^{in}_\ul(z) \} - \sum_\al [\la_\al(z), x_\al(z)]\right) .  
\end{align}

\subsection{Decomposition of $L\bar\G$}

Let us define $\G^{in}$ as the Lie subalgebra of $L\bar\G$
equal to $t\bar\G[[t]]$. Then 

\begin{lemma}
If $\Den(\uP,\ul)\neq 0$, we have a direct sum decomposition
$$
L\bar\G = \wt\G^{out}_{\ul} \oplus \G^{in} . 
$$
\end{lemma} 

{\em Proof.} $\wt\G^{out}_\ul$ is the space of all sums 
\begin{equation} \label{gestalt}
y(z) + \sum_i A_i (- e +  \ell_i h + \ell_i^2 f)
[ r^{(P_i)}(z) - r^{(P_0)}(z)] ,
\end{equation} 
where $z\mapsto y(z)$ is a regular map from $C \setminus \{P_0\}$ to $\bar\G$, 
and the $A_i$ are tangent vectors at $P_i$. 

The intersection $\wt\G^{out}_\ul \cap \bar\G[[t]]$ is therefore equal
to   the space of vectors of the form (\ref{gestalt}), such that  $y(z)
- [\sum_i A_i (- e + \ell_i h + \ell_i^2 f) ] r^{(P_0)}(z)$ is regular
at  $P_0$. Since there is no nonconstant function $\upsilon$ on the
universal cover $\wt C$ of $C$, regular everywhere and with
transformation properties $\upsilon\circ \gamma_{A_a} = \upsilon,
\upsilon\circ \gamma_{A_a} - \upsilon =$ constant, $y$ should be
constant and  $\sum_i A_i (- e + \ell_i h + \ell_i^2 f)$ should vanish.
Now the equations
$$
\sum_i A_i \omega_a(P_i) = \sum_i A_i \omega_a(P_i) \ell_i = \sum_i A_i \omega_a(P_i)
\ell_i^2 = 0
$$
imply that the $A_i$ are all zero. Therefore  $\wt\G^{out}_\ul \cap
\bar\G[[t]] = \bar\G$, so    $\wt\G^{out}_\ul \cap \G^{in} = 0$. 

Let us now prove that $\wt\G^{out}_{\ul} + \G^{in} = L\bar\G$. We have 
$$
L\bar\G / (\G^{in} + \G^{out}_\ul) = H^1(C, \End(\cE)_0(-P)) ,  
$$
and since $\cE$ is stable, we have $\dim H^1(C, \End(\cE)_0(-P)) = 3g$, so 
\begin{equation} \label{pidyon}
\dim L\bar\G / (\G^{in} + \G^{out}_\ul) = 3g.
\end{equation} 
On the other hand, we have $\dim\wt\G^{out}_\ul / \G^{out}_\ul = 3g$. 
Since the sum of $\G^{in}$ and $\wt\G^{out}_\ul$ is direct, we have also
 $\dim(\G^{in} + \wt\G^{out}_\ul) / (\G^{in} + \G^{out}_\ul) = 3g$.
Together with  (\ref{pidyon}), this shows $L\bar\G = \G^{in} +
\wt\G^{out}_\ul)$.   \hfill \qed\medskip

\subsection{The $\ell$-operators}

Let us set $\G_\Omega = \bar\G\otimes \CC((t))dt$.  There is a pairing
$\langle , \rangle_{L\bar\G\times\G_\Omega} :  L\bar\G\times \G_\Omega
\to\CC$, defined as the tensor product of $\langle, \rangle_{\bar\G}$ 
and the  residue pairing at $P_0$. 

Let us set  $\delta(z,w) dz = \sum_{i\in\ZZ} z^{i}w^{-i-1} dz$. Then the
canonical  element $L\bar\G\hat\otimes\G_\Omega$ is sent by the
isomorphism $\G_\Omega  \hat\otimes L\bar\G \to \bar\G^{\otimes
2}\otimes \CC[[z^{\pm1},w^{\pm1}]]dz$ to $\ell^{tot}(z,w) dz =  \Cas
\delta(z,w) dz$.  

We define $\ell^{in}_{\ul}(z,w)dz$ and $\ell^{out}_{\ul}(z,w)dz$ as the
images  in $\G_\Omega \otimes L\bar\G$ of the canonical elements of 
$(\wt\G_{\ul}^{out})^{\perp}\hat\otimes \G^{in}$ and 
$(\G^{in})^{\perp}\hat\otimes \wt\G_{\ul}^{out}$. 

\begin{prop}
Let us set   
\begin{align*}
& \ell^{rat}_{\ul}(z,w)dz = \Cas G(z,w) dz  
\\ & 
- {1\over {\Den(\uP,\ul})} \sum_{i = 1}^{3g} 
[- \pa^2_{\ell_i}\Den(\uP,\ul)_{| \ell_i = 0} e
+ {1\over 2}\pa_{\ell_i}\Den(\uP,\ul)_{| \ell_i = 0} h  
+ {1\over 2}\Den(\uP,\ul)_{| \ell_i = 0} f ]_{|P_i = z}
\\ & 
\otimes  ( - e + \ell_i h + \ell_i^2 f)G(P_i,w) dP_i  
\end{align*}
If $f(z,w)dz$ is a rational form on 
$C\times C$, let us denote by $f(z,w)dz_{|z\ll w}$  and $f(z,w)dz_{|w\ll z}$ 
its formal expansions for $z,w$ in the neighborhood of $P_0$ with 
$z\ll w$, resp.\ $w\ll z$. Then we have
$$
\ell^{out}_{\ul}(z,w)dz = - \ell^{rat}_{\ul}(z,w)dz_{|z\ll w}, \quad 
\ell^{in}_{\ul}(z,w)dz = \ell^{rat}_{\ul}(z,w)dz_{|w\ll z} . 
$$
\end{prop}

\subsection{Expression of $T(z)$ in terms of $\ell$-operators}

We will denote by $\iota$ the embedding of $[\bar\G\otimes
\Omega_z]\hat\otimes [\bar\G\otimes \CC((w))]$ in $\bar\G\otimes
[\G\hat\otimes \Omega_z]$ by sending  $[a\otimes \omega(z)]\otimes
[b\otimes f(w)]$ to  $a\otimes [(b\otimes f,0)\otimes \omega]$. (So
$\iota$ preserves the order in the tensor product $\bar\G\otimes\bar\G$,
but exchanges variables $z$ and $w$.)  

Let us denote by $\ell^{in}_{\ul}(z)$ and $\ell^{out}_{\ul}(z)$ the images of 
$\ell^{in}_{\ul}(z,w)dz$ and $\ell^{out}_{\ul}(z,w)dz$ by $\iota$. 
We have $\ell^{in}_{\ul}(z) +  \ell^{out}_{\ul}(z) = \ell(z)$. 

Let us define $\la(z)$ and the image by $\iota$ of
$\ell^{rat}_{\ul}(z,w)dz - \Cas G(z,w)dz$ and  define $\la_\al(z)$ by
$\la(z) = \sum_\al \la_\al(z) \otimes x_\al$.  Write also $\la(z) =
\sum_{u,v,\omega,\phi} u[\phi] \omega(z) \otimes v$.  Then 
$$ 
\sum_\al [\la_\al(z), x_\al(z)] = \sum \omega(z)\phi(z)
[u,v](z) + K \sum \omega(z) d\phi(z) \langle u,v \rangle_{\bar\G}.
$$ 
Therefore 
\begin{align*}
& \sum_\al [\la_\al(z), x_\al(z)] 
\\ & = - {1\over{\Den(\uP,\ul})}
\{ 
e(z) \sum_i [ - \pa_{\ell_i} \Den(\uP,\ul)_{|\ell_i = 0} 
+ 2 \ell_i \pa^2_{\ell_i} \Den(\uP,\ul)_{|\ell_i = 0}
]_{P_i = z} G(P_i,z) dP_i
\\ & 
+ h(z) \sum_i [ {1\over 2} \Den(\uP,\ul)_{|\ell_i = 0} 
- \ell_i^2 \pa^2_{\ell_i} \Den(\uP,\ul)_{|\ell_i = 0}
]_{P_i = z} G(P_i,z) dP_i
\\ & 
+ f(z) \sum_i [\ell_i \Den(\uP,\ul)_{|\ell_i = 0} 
- \ell_i^2 \pa_{\ell_i}\Den(\uP,\ul)_{|\ell_i = 0}]_{P_i = z}
G(P_i,z) dP_i \} 
\\ & 
- { k \over {\Den(\uP,\ul)} }
\sum_i d_z[G(P_i,z)dP_i] 
\\ &  \{ - {1\over 2} \Den(\uP,\ul)_{|\ell_i = 0}
+ \ell_i \pa_{\ell_i} \Den(\uP,\ul)_{|\ell_i = 0}
- \ell_i^2 \pa^2_{\ell_i} \Den(\uP,\ul)_{|\ell_i = 0} \}_{P_i = z} .  
\end{align*}

(\ref{new:sug}) then yields the following expression of the Sugawara
tensor: 

\begin{lemma} 
We have 
\begin{equation} \label{dmi}
T(z) ={1\over{2\kappa}} [(\tr\circ\rho) \otimes \id] 
\left( \ell^{out}_{\ul}(z)\ell(z) + \ell(z)\ell^{in}_{\ul}(z) 
+ a_\uP(\ul,z) \ell(z) \right) + {1\over 2\kappa} s_\uP(\ul,z) ,   
\end{equation}
where $a_\uP(\ul,z)$ is the element of $\bar\G\otimes 
\CC[\ell_i,\pa_{\ell_i}, i = 1,\ldots,3g]
[\Den(\uP,\ul)^{-1}]\otimes \Omega(C) $ equal to 
\begin{align*}
& a_\uP(\ul,z) = - {1\over{\Den(\uP,\ul)}}\sum_i G(P_i,z) dP_i 
\{ ({1\over 2} h + \ell_i f ) \Den(\uP,\ul)_{|\ell_i = 0} 
\\ &  
- (e + \ell_i^2 f) \pa_{\ell_i}\Den(\uP,\ul)_{|\ell_i = 0} 
+ ( 2 \ell_i e - \ell_i^2 h ) \pa^2_{\ell_i}\Den(\uP,\ul)_{|\ell_i = 0} 
\}_{P_i = z} , 
\end{align*}
and $s_\uP(\ul,z)$ is the element of $\CC[\ell_i,\pa_{\ell_i},i = 1,\ldots,3g]
[\Den(\uP,\ul)^{-1}] \otimes
\Omega^2(C)$ ($\Omega^2(C)$ is the space of rational sections of $\Omega_C^{\otimes 2}$) 
given by 
\begin{align*}
s_\uP(\ul,z) = & {k\over{\Den(\uP,\ul)}}
\sum_{i = 1}^{3g} d_z[G(P_i,z)dP_i]
\\ & \{ - {1\over 2} \Den(\uP,\ul)_{|\ell_i = 0}
+ \ell_i \pa_{\ell_i} \Den(\uP,\ul)_{|\ell_i = 0}
- \ell_i^2 \pa^2_{\ell_i} \Den(\uP,\ul)_{|\ell_i = 0} \}_{P_i = z}.  
\end{align*}
\end{lemma}

\subsection{}

For $x = x_e \otimes e + x_h \otimes h + x_f \otimes f$ an element of
$\wt\G^{out}_{(\ell_i)}$, let us set 
$$
\nu(x) = \sum_i (x_h + \ell_i x_e)(P_i).   
$$
We have then 
\begin{lemma}
For $v'$ any vector $\WW$,  and $x$ in  $\wt\G^{out}_{(\ell_i)}$, we have  $$
\langle \psi, \rho_\WW^{group}(\wt g_{\ul,\uP}) \rho_\WW[(x,0)] v' \rangle  
=  [\sum_{i = 1}^{3g} \mu_i(x) \pa_{\ell_i} 
- k \nu(x)] \langle \psi, \rho_\WW^{group}(\wt g_{\ul,\uP}) v' \rangle . 
$$
\end{lemma}

Recall that $\ell^{out}_{\ul}(z,w)dz$ belongs to $\G_\Omega \hat\otimes L\bar\G$. 
The expressions 
$$
\mu_{i,\uP}(\ul,z) := (\id\otimes \mu_i)(\ell^{out}_\ul(z,w)dz)
\ \on{and}\ \nu_\uP(\ul,z) := (\id\otimes \nu)(\ell^{out}_{\ul}(z,w)dz)
$$ 
make sense and belong to $\bar\G\otimes 
\CC[\ell_i, i  = 1,\ldots, 3g][\Den(\uP,\ul)^{-1}] \otimes \Omega(C)$. We find 
\begin{align*}
& \mu_{i,\uP}(\ul,z) =  (e - \ell_i h - \ell_i^2 f) G(z,P_i) dz  
 + {1\over{\Den(\uP,\ul)}} \cdot \\ & 
 \cdot \sum_{j\neq i} \ell_{ij}^2
[ \pa_{\ell_j}^2\Den(\uP,\ul)_{|\ell_j = 0} e
- {1\over 2} \pa_{\ell_j}\Den(\uP,\ul)_{|\ell_j = 0} h
- {1\over 2} \Den(\uP,\ul)_{|\ell_j = 0} f ]_{P_j = z}
G(P_j,P_i) dP_j 
\end{align*}
and 
\begin{align*}
&  \nu_\uP(\ul,z) = - \sum_i ({1\over 2} h + \ell_i f) G(z,P_i)dz 
+ {1\over{\Den(\uP,\ul)}} \cdot  
\\ & 
\cdot 
\sum_{i,j; i\neq j} \ell_{ji}
[ - \pa^2_{\ell_j}\Den(\uP,\ul)_{|\ell_j = 0} e  
+ {1\over 2} \pa_{\ell_j}\Den(\uP,\ul)_{|\ell_j = 0} h  
+ {1\over 2} \Den(\uP,\ul)_{|\ell_j = 0} f 
]_{P_j = z} G(P_j,P_i) dP_j .   
\end{align*}

\begin{lemma} \label{diff:expr}
Let us set 
$$
\ell^{diff}_\uP(z) := \sum_i \mu_{i,\uP}(\ul,z) \pa_{\ell_i} - k \nu_\uP(\ul,z) . 
$$
$\ell^{diff}_\uP(z)$ is an element of  $\bar\G\otimes
\CC[\ell_i,\pa_{\ell_i}, i = 1,\ldots,3g] [\Den(\uP,\ul)^{-1}]\otimes\Omega(C)$, 
and we have, for any vector $v'$ of $\WW$,  
$$
\langle \psi, \rho_\VV^{group}(\wt g_{\ul,\uP})
(\id\otimes \rho_\VV)(\ell^{out}_\ul(z)) v' \rangle 
= \ell^{diff}_\uP(z) \{   
\langle \psi, \rho_\VV^{group}(\wt g_{\ul,\uP})
v' \rangle \} . 
$$
\end{lemma}

Using the expression (\ref{dmi}) for $T(z)$, and the fact that 
$(\id\otimes \rho_\WW)(\ell^{in}_\ul(z)) (v)  = 0$, we find  
\begin{align*}
& 2\kappa \langle \psi , \rho_\WW^{group}(\wt g_{\ul,\uP}) \rho_\WW(T(z)) v_{0,k} \rangle
\\ & = 
\langle \psi , \rho_\WW^{group}(\wt g_{\ul,\uP}) \rho_\WW[ \sum_\al\ell^{out}_{\ul,
x_\al}(z) \ell_{y_\al}(z) + 
\sum_\al a_{x_\al,\uP}(\ul,z)\ell_{y_\al}(z)] v \rangle + s_\uP(\ul,z) f_\psi(\ul|v) , 
\end{align*}
where we set $\ell_{\ul,x}(z) = \langle x\otimes \id\otimes \id,
\ell_{\ul}(z)\rangle_{\bar\G}$, $\ell^{out}_{\ul,x}(z) = \langle x\otimes
\id\otimes \id, \ell^{out}_{\ul}(z) \rangle_{\bar\G}$
$a_{x,\uP}(\ul,z) = \langle x\otimes \id, a_\uP(\ul,z), \rangle_{\bar\G}$, 
for $x$ in $\bar\G$. 

Let us set, for $x$ in $\bar\G$, 
$\ell^{diff}_{\uP,x}(z) = \langle x\otimes \id\otimes \id,
\ell^{diff}_\uP(z)\rangle_{\bar\G}$. $\ell^{diff}_{x,_uP}(z)$ belongs then to 
$\CC[\ell_i,\pa_{\ell_i}, i = 1, \ldots, 3g][\Den(\uP,\ul)^{-1}]\otimes \Omega(C)$. 
Lemma \ref{diff:expr} then implies that 
\begin{align*}
& 2\kappa \langle \psi , \rho_\WW^{group}(\wt g_{\ul,\uP}) \rho_\WW(T(z)) v \rangle
\\ & = 
\sum_\al [\ell^{diff}_{\uP,x_\al}(z) + a_{\uP,x_\al}(\ul,z)] \{ 
\langle \psi , \rho_\WW[\ell^{out}_{\ul}(z)] v \rangle \}  
+ s_\uP(\ul,z) f_\psi(\ul,\uP|v) 
\\ &  = 
\{ [\sum_\al \ell^{diff}_{\uP,x_\al}(z)\ell^{diff}_{\uP,y_\al}(z) 
+ \sum_\al a_{\uP,x_\al}(\ul,z)\ell^{diff}_{\uP,y_\al}(z) + s_\uP(\ul,z) ]
f_\psi \}(\ul,\uP|v) . 
\end{align*}
We have therefore 
\begin{thm}
Let us set 
\begin{equation} \label{formula:T:diff}
T^{diff}_\uP(z) = 
\sum_\al \ell^{diff}_{\uP,x_\al}(z)\ell^{diff}_{\uP,y_\al}(z) 
+ \sum_\al a_{\uP,x_\al}(\ul,z)\ell^{diff}_{\uP,y_\al}(z) + s_\uP(\ul,z) . 
\end{equation}
Then $T^{diff}_\uP(z)$ belongs to   $\CC[\ell_i,\pa_{\ell_i}, i = 1, \ldots, 3g]
[\Den(\uP,\ul)^{-1}]\otimes H^0(C,\Omega^2_C(2\sum_{i=1}^{3g} P_i))$. In other words, 
it is a quadratic differential in $z$, regular on $C$ except for double poles at 
the $P_i, i = 1,\ldots, 3g$, with coefficients differential operators in $\ul$. 

There are differential operators $T^{diff}_{\uP,i}$ in $\ul$, 
polynomial in $k$ of degree $\leq 1$, such that the 
expansion of  $T^{diff}_\uP(z)$ for $z$ near $P_i$  is 
\begin{equation} \label{expansion:T:diff}
T^{diff}_\uP(z) = 2\kappa k ({1\over 2} {{dz_{P_i}}\over{z_{P_i}}})^2 
+ \kappa T^{diff}_{\uP,i} {{(dz_{P_i})^2}\over{z_{P_i}}} + O(1) .  
\end{equation}

For any $\G^{out,ext}$-invariant form $\psi$ on $\WW$, 
we have 
\begin{equation} \label{action:T}
\langle \psi, \rho_{\WW}^{group}(\wt g_{\ul,\uP}) \rho_\WW(T(z)) 
v \rangle = {1\over{2\kappa}}(T^{diff}_\uP(z)f_\psi)(\ul|v) . 
\end{equation}
\end{thm}

{\em Proof of Thm.} (\ref{action:T}) follows from the above computations. 
Let us now show the properties of $T^{diff}_\uP(z)$.   
We will use the following Lemma

\begin{lemma} \label{oberland}
Let us set 
\begin{align} \label{diff:T:wtT}
& \wt T_\uP^{diff}(z) = T_\uP^{diff}(z) +  \kappa 
{{db_\uP}\over{b_\uP}}(z) \ell^{diff}_{\uP,h}(z) \\ & \nonumber + 
2\kappa [\sum_i \ell_i^{-1} \{ d(d_{i,\uP}) - d_{i,\uP} {{db_\uP}\over
{b_\uP}}\}(z)] \ell^{diff}_{\uP,e}(z) + 2k\kappa ({1\over 2}
{{db_\uP}\over{b_\uP}})^2(z) . 
\end{align}
Then we have 
$$
2\kappa \langle \psi, \rho_\WW(T(z)) \rho_\WW^{group} (\wt g_{\ul,\uP}) v \rangle  
=  (\wt T_\uP^{diff}(z) f_\psi)(\ul,\uP|v) . 
$$
\end{lemma}

Let $k$ be an arbitrary complex number and  let $\wt W$ be any 
$\bar\G$-module.   Define $\wt\WW$ as the Weyl module over $U\G$ equal
to $U\G\otimes_{U\G_+} \wt W$ (see sect.\ \ref{sect:def:weyl}). Let us
denote by $\rho_{\wt\WW}$ the corresponding  algebra map from $U\G$ to
$\End(\wt\WW)$. For $v$ in $\wt W$, we denote also by $v$ the vector
$1\otimes v$ of $\wt\WW$. 

Let us fix $\ul_0 = (\ell_{i,0})_{i = 1,\ldots, 3g}$ 
in $(\CC^\times)^{3g}$, such that $\Den(\uP,\ul_0) \neq 0$. 
Consider $\ul$ as a formal variable near $\ul_0$ and let us define 
$\gamma_{\ul_0,\ul,\uP}$ as the element of  
$U\G[[\ul - \ul_0]]$ equal to 
$$
\gamma_{\ul_0,\ul,\uP} = \exp(e[\sum_i (\ell_{i,0}^{-1} - \ell_i ^{-1}) d_{i,\uP}])
z(\prod_i (\ell_i / \ell_{i,0})^k) . 
$$ 
When $k$ is an integer $\geq 0$, we have 
$\wt g_{\uP,\ul_0}^{-1}
\wt g_{\uP,\ul} = \gamma_{\ul_0,\ul,\uP}$. 
Then the analogue of Lemma \ref{oberland} and of formula (\ref{action:T}) 
for the Weyl module $\wt\WW$ are:  

\begin{lemma}  Let $\psi_0$ be a
$\Ad(g_{\uP,\ul_0})(\G^{out})$-invariant form on $\wt\WW$. For any $v'$
in $\wt\WW$,  $\langle \psi_0, \rho_{\wt\WW}[\Ad( g_{\uP,\ul_0}^{-1} )(  T(z))
\gamma_{\ul_0,\ul,\uP}]( v') \rangle$ and $\langle \psi_0, 
\rho_{\wt\WW}[\gamma_{\ul_0,\ul,\uP}]( v') \rangle$ are formal functions in $\ul -
\ul_0$, and we have, for any $v$ in $\wt W$,    
\begin{equation} \label{generic:id}
\langle \psi_0, \rho_{\wt\WW}[\Ad( g_{\uP,\ul_0}^{-1} )( 
T(z)) \gamma_{\ul_0,\ul,\uP}]( v)
 \rangle = \wt T_\uP^{diff}(z)\{\langle \psi_0,  \rho_{\wt\WW}[\gamma_{\ul_0,\ul,\uP}](v)
\rangle \}
\end{equation} 
and 
\begin{equation} \label{generic:id'}
\langle \psi_0,  \rho_{\wt\WW}[\gamma_{\ul_0,\ul,\uP}  T(z)]( v) \rangle 
= T_\uP^{diff}(z)
\{\langle \psi_0,  \rho_{\wt\WW}[\gamma_{\ul_0,\ul,\uP}]( v) \rangle \}
\end{equation} 
(equalities in $\CC[[\ul - \ul_0]]$).  
\end{lemma}

The interest of identities (\ref{generic:id}) and (\ref{generic:id'}) 
lies in the fact that any formal function of $\ul - \ul_0$ is of
the form  $\langle \psi_0,  \rho_{\wt\WW}[\gamma_{\ul_0,\ul,\uP}]( v)
\rangle$. More precisely, we have: 

\begin{lemma} \label{chapel:hill}
The correlation functions map $\psi_0 \mapsto (v\mapsto 
\langle \psi_0, 
\rho_{\wt\WW}[\gamma_{\ul_0,\ul,\uP}](v) \rangle) $ defines an isomorphism 
from the space of $\Ad(g_{\uP,\ul_0})(\G^{out})$-invariant forms on $\wt\WW$
to $\Hom_{\bar\G}(\wt\WW, \CC[[\ul - \ul_0]])$. 
\end{lemma}

Here $\CC[[\ul - \ul_0]]$ is endowed with the action of $\bar\G$ defined
 by $e\mapsto \sum_i{{\pa}\over{\pa\ell_i}}$, $h\mapsto   \sum_i -
(2\ell_i {{\pa}\over{\pa\ell_i}} + k)$,  $f\mapsto \sum_i - (\ell_i^2
{{\pa}\over{\pa\ell_i}} + k \ell_i)$.

Lemma \ref{chapel:hill} follows from the fact that the map 
$$
U\G^{out}_\ul \otimes \CC[e[d_{i,\uP}]] \otimes U(\bar\G[[t]])
\to UL\G, 
$$
induced by the multiplication, is surjective, and its kernel is 
generated by the relations
$$
e[1] - \sum_i \ell_i^2 e[d_{i,\uP}] \in \G^{out}_\ul, 
h[1] - \sum_i 2\ell_i e[d_{i,\uP}] \in \G^{out}_\ul, 
f[1] + \sum_i e[d_{i,\uP}] \in \G^{out}_\ul. 
$$
 
Taking for example $\wt W = \CC[[\ul - \ul_0]]$, we see that any 
element of $\CC[[\ul - \ul_0]]$ can be obtained as a
$\langle \psi_0,  \rho_{\wt\WW}[\gamma_{\ul_0,\ul,\uP}](v) \rangle$. 

For any vector field $\xi$ regular on $C\setminus \{P_0\}$, $T[\xi]
= \res_{P_0}(\xi(z)T(z))$ belongs to $\G^{out} U\G$ (see \cite{EF:sols}). 
It follows that the right side of (\ref{generic:id}) 
is a quadratic differential in $z$, regular outside $P_0$. 

It follows that  $\wt T^{diff}_\uP(z)$ is a quadratic differential in
$z$,  regular outside $P_0$, with coefficients differential operators in
$\ul$.  Moreover, formula (\ref{formula:T:diff})  expressing
$T^{diff}(z)$ implies that $T^{diff}(z)$ is regular at $z = P_0$. 

The other functional properties of $T^{diff}(z)$ follow from 
formula (\ref{diff:T:wtT}) expressing $T^{diff}(z)$
in term of $\wt T^{diff}(z)$: 

-- $d\log b_\uP$ and $\ell^{diff}_\uP$ are univalued differentials in $z$, therefore 
$T^{diff}(z)$ is a univalued quadratic differential; 

-- the expansions 
$$
\ell^{diff}_{\uP,h}(z) = 2 {{dz_{P_i}}\over{z_{P_i}}} (-\ell_i \pa_{\ell_i}
+ {k\over 2}) + O(1) , \quad  
\ell^{diff}_{\uP,e}(z) =  {{dz_{P_i}}\over{z_{P_i}}} (-\ell_i^2 \pa_{\ell_i}
+ k\ell_i) + O(1) ,   
$$
and $s_\uP(z) = ({{dz_{P_i}}\over{z_{P_i}}})^2 + O(z_{P_i}^{-1})$
for $z$ near $P_i$ imply that $T^{diff}(z)$ has the expansion 
$$
T^{diff}(z) = {1\over 2} \kappa k ({{dz_{P_i}}\over{z_{P_i}}})^2 + O(z_{P_i}^{-1}) 
$$
for $z$ near $P_i$. Expansions (\ref{expansion:T:diff}) follow. Since
the  coefficients of ${{(dz_{P_i})^2}\over{ z_{P_i}}}$ in these
expansions  are the same as those of $(T^{diff} - \wt T^{diff})(z)$, 
which are given by formula (\ref{diff:T:wtT}), they are polynomial in 
$k$, of degree $\leq 1$. 
\hfill \qed\medskip

{\em Proof of Lemma \ref{oberland}.}
Let us denote by $(U\G)_k$ the quotient of $U\G$ by the ideal
generated by $K - k\cdot 1$, and by $(U\G)_{k,compl}$ the completion
of $(U\G)_k$ w.r.t.\ the topology defined by the right ideals generated
by the $x[z^l],l\geq N$ (the Laurent coefficients of $T(z)$ belong 
to $(U\G)_{k,compl}$.) 

Then $\wh{G(\cK)}$ acts on $(U\G)_{k,compl}$ by the adjoint action; this action 
factors through $G(\cK)$; so  we have a sequence of group morphisms
$$
\wh{G(\cK)} \to G(\cK) \to \Aut_{algebra} [ (U\G)_{k,compl}] . 
$$
We denote by $\Ad$ any of the maps $G(\cK) \to \Aut_\CC [ (U\G)_{k,compl}]$
and $\wh{G(\cK)} \to \Aut_\CC [ (U\G)_{k,compl}]$. 

Moreover, it follows from the formula
$$
[T(z) , x[f]] = - (df)(z) x(z) 
$$
for $x$ in $\bar\G$, that if $f$ belongs to $\cK$, 
$$
\Ad( \pmatrix 1 & f \\ 0 & 1 \endpmatrix ) [T(z)]
= T(z) + (df)(z) e(z)  ,   
$$
and that for $\phi$ in $\CC[[z_{P_0}]]^\times$,  
$$
\Ad( \pmatrix \phi & 0 \\ 0 & \phi^{-1} \endpmatrix ) [T(z)]
= T(z) + {{d\phi}\over\phi}(z) h(z) + K ({{d\phi}\over\phi})^2(z) .    
$$
Moreover, we have also 
$$
\Ad(T_{\wh{G(\cK)}})[T(z)] = T(z) + {1\over 2}{{dz_{P_0}}\over
{z_{P_0}}}(z) h(z) + K ({1\over 2}{{d z_{P_0}}\over {z_{P_0}}})^2(z) . 
$$
Let then $\psi$ belong to $\cK^\times$; let us denote $v(\psi)$ the valuation of 
$\psi$ at $P_0$, and set $t[\psi] = T_{\wh{G(\cK)}}^{v(\psi)} \pmatrix \psi_0^{1/2} & 0 
\\ 0  & (\psi_0^{1/2})^{-1} \endpmatrix$, where $\psi_0 = z_{P_0}^{-v(\psi)}
\psi$ and $\psi_0^{1/2}$ is any square root of $\psi_0$, we have therefore
$$
\Ad(t[\psi])(T(z)) = T(z) + {1\over 2} {{d\psi}\over{\psi}}(z) h(z)
+ K ({1\over 2}{{d\psi}\over{\psi}})^2(z).  
$$

We have therefore 
\begin{align*}
& \Ad(\wt g_{\ul,\uP})^{-1}(T(z))  = 
\Ad(\iota( \pmatrix 1 & \sum_i \ell_i^{-1} d_{i,\uP} 
\\ 0 & 1 \endpmatrix ) )[T(z) + {1\over 2}{{db_\uP}\over{b_\uP}}(z) h(z) 
+ K  ({1\over 2} {{db_\uP}\over {b_\uP}})^2(z) ] 
\\ & = 
T(z) + [\sum_i \ell_i^{-1} \{ d(d_{i,\uP}) - d_{i,\uP} {{db_\uP}\over {b_\uP}}\}(z)] e(z)  
+ {1\over 2}{{db_\uP}\over {b_\uP}}(z) h(z) + K ({1\over 2} {{db_\uP}\over {b_\uP}})^2(z).  
\end{align*}

So 
\begin{align*}
& 2\kappa \langle \psi, \rho_\WW(T(z)) \rho_\WW^{group} (\wt g_{\ul,\uP}) v \rangle  
\\ &  
= 2\kappa \langle \psi, \rho_\WW^{group} (\wt g_{\ul,\uP}) \rho_\WW(T(z)) v  \rangle  
\\ & 
+  2\kappa [\sum_i \ell_i^{-1} \{ d(d_{i,\uP}) - d_{i,\uP} {{db_\uP}\over {b_\uP}}\}(z)]
\langle \psi, \rho_\WW^{group} (\wt g_{\ul,\uP}) \rho_\WW(e(z)) v \rangle 
\\ & 
+ 2\kappa {1\over 2}{{db_\uP}\over b_\uP}(z)
\langle \psi, \rho_\WW^{group} (\wt g_{\ul,\uP}) \rho_\WW(h(z)) v \rangle 
\\ & 
+ 2\kappa k ({1\over 2} {{db_\uP}\over {b_\uP}})^2(z) 
\langle \psi, \rho_\WW^{group} (\wt g_{\ul,\uP}) v \rangle  . 
\end{align*}

On the other hand, we have for any $x$ in $\bar\G$, 
$$
\langle \psi, \rho_\WW^{group}(\wt g_{\ul,\uP}) \rho_\WW(x(z)) v \rangle
= (\ell^{diff}_{\uP,x}(z) f_\psi)(\ul,\uP|v) . 
$$
Lemma \ref{oberland} follows. 
\hfill \qed\medskip

\section{Explicit form of the Beilinson-Drinfeld operators}
\label{sect:bd}

\subsection{Pull-back of the $\det$ line bundle}

Recall that the line bundle $\det$ is defined over 
$\Bun_2(C)_{\cO(\sum_i P_i)}$ as $\det R\Gamma$ (\cite{Beauville}).

{\em Proof of Prop.\ \ref{prop:image:det}.}
A line bundle over $(\CC P^1)^N$ is necessarily of the 
form $\cO(k_1)\boxtimes\cdots\boxtimes \cO(k_N)$. To determine  the
$k_i$, let us compute the restriction of  $\cH_\uP^*(\det)$  to
$(\ell_1,\ldots,\ell_{i-1})\times \CC P^1\times  (\ell_{i+1},\ldots,
\ell_N)$, for generic $\ell_i$. For $\ul$ generic, the fiber of  $\cH_\uP^*(\det)$ at $\ul$
is the determinant of  
$\{(\al_i)\in  \CC^{3g} | \sum_i \al_i \omega_a(P_i) 
= \sum_i \al_i \ell_i \omega_a(P_i) = 0 \}$, which we view as 
$$
\{(\al_j,\al'_i)_{j\neq i}\in  \CC^{3g} | \sum_{j\neq i} \al_j \omega_a(P_j)  
+ \al'_i \ell_i^{-1}\omega_a(P_i) 
= \sum_{j\neq i} \al_j \ell_j \omega_a(P_j) 
+ \al'_i \omega_a(P_i)
 = 0 \}
$$ 
in the neighborhood of $\ell_i = \infty$. It follows that 
the restriction of  $\cH_\uP^*(\det)$  to
$(\ell_1,\ldots,\ell_{i-1})\times \CC P^1\times  (\ell_{i+1},\ldots,
\ell_N)$ is isomorphic to $\cO(1)$, therefore all 
$k_i$ are equal to $1$. 
\hfill \qed\medskip

\subsection{Explicit form of the Beilinson-Drinfeld operators}

The Beilinson-Drinfeld operators are differential operators 
acting on sections of 
$\det^{-2}$. They are defined as the action on  $\Bun_2(C)_{\cO(\sum_i
P_i)}$ of the central elements of the critical level quotient $U\G / (K
+ 2\cdot 1)$ of the enveloping algebra of $\G$. A generating series for
these central  elements is the Sugawara series $T(z) (dz)^2= \sum_n T_n
z^{-n-2}(dz)^2$;  let $T_n^{(BD)}$ be the action of $T_n$ by
differential operators.  Then the series $\sum_n T_n^{(BD)}
z^{-n-2}(dz)^2$ is a linear  combination of the quadratic differentials
on $C$.  Let $(\omega_\al^{(2)})_{a = 1,\ldots, 3g - 3}$  be a basis of
$H^0(C,\Omega_C^{\otimes 2})$ and let us set  $\sum_n T_n^{(BD)}
z^{-n-2}(dz)^2  =  \sum_\al T_\al^{BD} \omega^{(2)}_\al$. Then the
$T_n^{(BD)}$ are all linear combinations of the $T_\al^{BD}$.

{\em Proof of Thm.\ \ref{thm:bd}}
The expansions (\ref{expansion:T:diff}) imply that 
$T^{diff}(z)$ has no poles at the $P_i$. The fact that the operators 
$T^{diff}_{\uP,\al}$ commute follows from formula (\ref{generic:id'}),
the fact that any 
formal series in $\CC[[\ul - \ul_0]]$ can be obtained as a correlation
function $\langle \psi_0, \rho[\gamma_{\ul_0,\ul,\uP}](v) \rangle$, and the
commutativity of the components of the Sugawara tensor when $\kappa = 0$. 
\hfill \qed\medskip

\section{Explicit form of the KZB connection} \label{sect:kzb}

\subsection{Dependence of $i_\uP(\psi)$ in $\uP$}

Let us denote by $(C^{3g})_{P_0}$ the subspace of 
$(C \setminus \{P_0\})^{3g}$, formed of the systems $\uP = 
(P_i)_{i = 1,\ldots, 3g}$, such that  $\sum_i A(P_i) = 3g A(P_0)$.  

Let us denote by $i_{\uP}$ the linear map 
$$
i_{\uP} : 
(\WW^*)^{\G^{out}} \to \Hom_{\bar\G}(W, \otimes_{i=1}^{3g} \CC[\ell_i]_{\leq k}) , 
$$ 
indexed by $\uP$ in $(C^{3g})_{P_0}$: we have 
$i_{\uP}(\psi) = [v\mapsto f_{\psi}(\ul,\uP|v)]$. 
In this section, we describe the dependence of $i_{\uP}(\psi)$
on $\uP$, when $\psi$ is fixed. For this, we first compare the 
maps $\cH_\uP$ and $\cH_{\uP'}$, for $\uP'$ infinitely close to $\uP$. 

\subsubsection{Infinitesimal comparison of $\cH_\uP$ and $\cH_{\uP'}$}

Let us fix variations $\delta P_i$ of the points $P_i$, such
that $\delta(\sum_i A(P_i)) = 0$. Let us fix $\ul$ such that $\Den(\uP,\ul) \neq 0$. 
It follows from Thm.\ \ref{thm:local:isom} that one can find $\ul'$ infinitely 
close to $\ul$, such that $\cH_{\uP}(\ul) =  \cH_{\uP'}(\ul')$. In this
section, we compute the variation $\delta\ul = \ul' - \ul$.  

Let us denote by 
$M(z)$ the matrix corresponding to the isomorphism from 
$\cH(\uP',\ul')$ to $\cH(\uP,\ul)$. We have 
$$
M(z) = \id_{\CC^2} + \sum_i 
\pmatrix 1 \\ \ell_i\endpmatrix  \pmatrix \al_i & \beta_i
\endpmatrix  r^{(P_i)}(z) , 
$$
where $\al_i = \al_i^{(1)} + \al_i^{(2)} + \ldots$ and $\beta_i = \beta_i^{(1)} + 
\beta_i^{(2)} + \ldots$ are infinitesimal vectors at $P_i$ ($\al_i^{(1)}$ is 
the linear in $\delta P_i$ part of $\al_i$, $\al_i^{(2)}$ the quadratic part, etc.), 
satisfying the  equations
\begin{equation} \label{smth}
\sum_i \omega_a(P_i) \pmatrix 1 \\ \ell_i \endpmatrix 
\pmatrix \al_i & \beta_i \endpmatrix = 0  
\end{equation}
for each $a$, and $M(P_i + \delta P_i) \pmatrix 1 \\ \ell_i + \delta \ell_i \endpmatrix 
= 0$, for any $i$, which implies at first order   
$$
\pmatrix 1 \\ \ell_i + \delta \ell_i \endpmatrix 
+ \sum_{j\neq i} r^{(P_j)}(P_i) (\al_j^{(1)} + \ell_i \beta_j^{(1)}) 
\pmatrix 1 \\ \ell_j \endpmatrix 
+ [r^{(P_i)}(P_i + \delta P_i) (\al_i^{(1)} + \ell_i\beta_i^{(1)}) + \kappa_i]
\pmatrix 1 \\ \ell_i \endpmatrix = 0 , 
$$ 
where we set 
$\kappa_i = r^{(P_i)}(P_i + \delta P_i) (\al_i^{(2)} + \ell_i \beta_i^{(2)})$. 

The main terms of the left side of this equation should cancel, 
therefore 
$$
\beta_i^{(1)}  = - {1\over{\ell_i}} (\al_i^{(1)} + \delta P_i).  
$$
(\ref{smth}) then yields the system 
\begin{equation} \label{system:dependence}
\sum_i \omega_a(P_i) \al_i^{(1)} = 0,   \quad
\sum_i \omega_a(P_i) \ell_i \al_i^{(1)} = 0, \quad  
\sum_i \omega_a(P_i) \ell_i^{-1} [\al_i^{(1)} + \delta P_i] = 0,  
\end{equation}
as well as the equation $\sum_i \omega_a(P_i)\delta P_i = 0$ 
for any $a = 1,\ldots, g$, which is satisfied because $\delta(\sum_i A(P_i)) = 0$. 
We then get 
$$
\kappa_i = - \sum_{j\neq i} r^{(P_j)}(P_i) (\al_j^{(1)} + \ell_i
\beta_j^{(1)}) 
$$ 
and 
$$
\delta \ell_i = - \sum_{j\neq i} r^{(P_j)}(P_i) \ell_j 
(\al_j^{(1)} + \ell_i \beta_j^{(1)}) - \kappa_i \ell_i
= 
\sum_{j\neq i} r^{(P_j)}(P_i) \ell_{ij} 
(\al_j^{(1)} + \ell_i \beta_j^{(1)}) . 
$$
(\ref{system:dependence}) yields 
$$
\al_i^{(1)} = {1\over{\Den(\uP,\ul)}} \sum_j - (\ell_i / \ell_j) \delta P_j 
\Den(\uP,\ul)_{|\ell_i = 0, P_i\to P_j} .  
$$
It follows that 
\begin{equation}
\delta \ell_i = {1\over{\Den(\uP,\ul)}} 
\sum_{j\neq i} r^{(P_j)}(P_i) \ell_{ij}^2
\sum_l {1\over{\ell_l}} \delta P_l
\Den(\uP,\ul)_{|\ell_j = 0, P_j\to P_l}  
- \sum_{j\neq i} r^{(P_j)}(P_i) \ell_{ij} {{\ell_i}\over{\ell_j}}
\delta P_j .  
\end{equation}

\subsubsection{Dependence of $i_\uP(\psi)$ on $\uP$}

Let us preserve the notation of the previous section. 

Let us set 
$$
m^{in}(z) = \sum_i \pmatrix \beta^{(1)}_i \\ -\al^{(1)}_i \endpmatrix 
\pmatrix \ell_i & -1\endpmatrix r^{(P_i)}(z). 
$$
The Laurent expansion of $m^{in}(z)$ at $P_0$ belongs to 
$\GL_2[[z_{P_0}]]$.  It follows from the previous section that 
there exists an element $m^{out}$ of $\GL_2\otimes
R$  such that 
$$
g_{\uP + \delta\uP,\ul + \delta\ul}  
= \exp(m^{out})   g_{\uP,\ul} \exp(m^{in}) . 
$$

Let us denote by $G(\cK)\rtimes \ZZ$ the quotient of $\wh{G(\cK)}
\rtimes \ZZ$ by its center,  and let $\wt g'_{\uP,\ul}$ be the image of 
$\wt g_{\uP,\ul}$ by the projection $\wh{G(\cK)}\rtimes \ZZ \to
G(\cK)\rtimes \ZZ$. Let $\pi$ be the projection map from $\GL_2$ to $\SL_2$
and set 
$$
m^{\prime out} = (\pi\otimes\id)(m^{out}), \quad  
m^{\prime in} = (\pi\otimes\id)(m^{in}). 
$$ 
Then we have 
$$
\wt g'_{\uP + \delta\uP,\ul + \delta\ul}  
= \exp(m^{\prime out})   \wt g'_{\uP,\ul} \exp(m^{\prime in}) ;  
$$
the Lie algebraic meaning of this equality is that 
\begin{equation} \label{prebalus}
\wt g'_{\uP + \delta\uP,\ul + \delta\ul}  
\wt g_{\uP,\ul}^{\prime -1} = m^{\prime out} + 
\wt g'_{\uP,\ul} m^{\prime in} \wt g_{\uP,\ul}^{\prime -1} . 
\end{equation}

Let us denote by $x^{in}$ and $x^{out}$ the images of $m^{\prime in}$ 
and  $m^{\prime out}$ by the map $L\bar\G\to\G$, $x\otimes f\mapsto x[f]
(= (x\otimes f,0))$. Then $x^{in}$ belongs to the linear span 
$\G^{+,no\ K}$ of the 
$x[t^i],i\geq 0$, $x\in\bar\G$, and $x^{out}$ belongs to $\G^{out,ext}$. 

Let us denote by $\kappa$ the map from $\G$ to $\CC$ defined by 
$\kappa(x\otimes f, aK) = ak$. 

\begin{prop}
Let us set 
$$
Z = \kappa \{ \wt g_{\uP + \delta\uP, \ul + \delta\ul} 
\wt g_{\uP, \ul}^{-1}  \}
- \kappa \{ \wt g_{\uP, \ul} x^{in}\wt g_{\uP, \ul}^{-1}\} . 
$$
Then we have 
\begin{equation} \label{Z}
Z = k \sum_{j\neq i} \al_i^{(1)} {{\ell_i}\over{\ell_j}} 
r^{(P_i)}(P_j) + k \sum_i {{\delta\ell_i}\over{\ell_i}}, 
\end{equation}
and 
\begin{equation} \label{balus}
\wt g_{\uP + \delta\uP,\ul + \delta\ul} 
\wt g_{\uP,\ul}^{-1}  = 
x^{out} + \wt g_{\uP,\ul} x^{in}\wt g_{\uP,\ul}^{-1} + Z, 
\end{equation}
therefore 
\begin{equation} \label{cobalus}
\wt g_{\uP + \delta\uP,\ul + \delta\ul}  
= \exp(x^{out})   \wt g_{\uP,\ul} \exp(x^{in})  z(\exp(Z)) , 
\end{equation}
up to quadratic terms in the $\delta\ell_i$.  
\end{prop}

{\em Proof.} (\ref{balus}) follows from (\ref{prebalus}), the 
definition of $Z$ and the fact that $\kappa(x^{out}) = 0$. 
 
Let us now compute $Z$. We have 
$ \kappa (\wt g_{\uP + \delta\uP, \ul + \delta\ul} 
\wt g_{\uP, \ul}^{-1})  = \kappa (t[b_{\uP + \delta\uP}^{-1/2}] 
t[b_{\uP}^{-1/2}]^{-1}) 
+ k \sum_i {{\delta\ell_i}\over{\ell_i}}
= \kappa (T_{\wh{G(\cK)}}^{3g} t[(b^{0}_{\uP + \delta\uP})^{-1/2} 
(b^{0}_{\uP})^{1/2}] T_{\wh{G(\cK)}}^{-3g}) 
+ k \sum_i {{\delta\ell_i}\over{\ell_i}}
= - {{3gk}\over 2} [{{b_{\uP + \delta\uP}}\over{b_\uP}}(P_0) - 1]
+ k \sum_i {{\delta\ell_i}\over{\ell_i}}$. 

On the other hand, let us set $x^{in} = e[x^{in}_e] 
+ h[x^{in}_h]  + + f[x^{in}_f]$. We have 
$$
x^{in}_f = - \sum_i \al^{(1)}_i \ell_i r^{(P_i)}(z), \quad
x^{in}_h = \sum_i 
{{\beta^{(1)}_i \ell_i - \al^{(1)}_i}\over 2} r^{(P_i)}(z).  
$$
We have then $\kappa(\wt g_{\uP,\ul} x^{in} \wt g_{\uP,\ul}^{-1})
= \kappa (t[b_\uP^{-1/2}] 
\{h[x^{in}_h]   + [e[ - \sum_i \ell_i^{-1} d_{i,\uP}], f[x^{in}_f]]\}
t[b_\uP^{-1/2}]^{-1})$. After some computation, we find 
$$
\kappa(\wt g_{\uP,\ul} x^{in} \wt g_{\uP,\ul}^{-1}) = 
k \sum_{j\neq i} \al_i^{(1)} {{\ell_i}\over{\ell_j}} r^{(P_i)}(P_j)
- {{3g k}\over 2} \sum_i \delta P_i r^{(P_i)}(P_0) . 
$$
Therefore 
$$
Z = - {{3gk}\over 2} [{{b_{\uP + \delta\uP}}\over{b_\uP}}(P_0) - 1]
+ k \sum_{j\neq i} \al_i^{(1)} {{\ell_i}\over{\ell_j}} r^{(P_i)}(P_j)
- {{3gk}\over 2} \sum_i \delta P_i r^{(P_i)}(P_0)  
+ k \sum_i {{\delta\ell_i}\over{\ell_i}} . 
$$
Since 
$$
{{\delta b_\uP}\over{b_\uP}}(z) + \sum_i \delta P_i r^{(P_i)}(z) 
$$
vanishes for any $z$, 
$$
[{{b_{\uP + \delta\uP}}\over{b_\uP}}(P_0) - 1]
+ \sum_i \delta P_i r^{(P_i)}(P_0)  = 0 ; 
$$
(\ref{Z}) follows.  
\hfill \qed\medskip 

Prop.\ \ref{variation:f:psi} follows from the insertion 
of (\ref{cobalus}) in correlation functions. \hfill \qed\medskip 

\section{Filtration of conformal blocks} \label{sect:filt}

Let again $\bar\G$ be an arbitrary semismple Lie algebra, and let 
$\VV$ be an integrable module over the affine Kac-Moody algebra 
$\G$.  Let $CB(\VV)$ be the space of conformal blocks of $\VV$.  
The maps $\cH_{\uP,\uchi}$ naturally lead to a filtration of 
$CB(\VV)$. It is defined as follows. 

It follows from \cite{BL} that conformal  blocks can be identified with
sections of line bundles over  $\Bun_G(C)$. Define   
$\Fil_{\uP,\uchi}$ as the
space of conformal blocks, vanishing on $\Imm \cH_{\uP,\uchi}$. 
We then have a inclusions $\Fil_{\uP,\uchi} \subset 
\Fil_{\uP,\uchi'}$ if $\chi_i \leq\chi'_i$ for any $i$. 

It is easy to define analogues $\Fil_{\uP,\uchi,P}$ of 
$\Fil_{\uP,\uchi}$, which are conformal blocks vanishing on 
modifications of the fixed (nontrivial) bundle $P$.  

Moreover, one may conjecture the following behavior of the 
$\Fil_{\uP,\uchi}$: 

1) the $\Fil_{\uP,\uchi}$ form a bundle over $\cM_{g,N}$
with a flat connection; and the $\Fil_{\uP,\uchi,P}$ 
form a bundle over the joint moduli space of triples $(C,\uP,P)$
of a curve, marked points and a principal $G$-bundle over it;

2) when $P_i$ and $P_j$ both tend to a point $Q$, the limit of 
$\Fil(\uP,\uchi)$ is $\Fil(\uP',\uchi')$, where 
$\uP' = ((P_\al)_{\al\neq i,j}, Q)$ and $\uchi' = 
((\chi_\al)_{\al\neq i,j}, \chi_i + \chi_j)$. 

It would follow from this that the dimension of 
$\Fil_{\uP,\uchi,P}$  only depends on the sum $\sum_i \chi_i$, 
so that there is a filtration of $CB(\VV)$ by $P_+$. It would be
interesting to see whether the formulas for the corresponding 
$q$-dimensions agree with those of \cite{Feigin:Loktev,Leclerc}.

\end{document}